\documentclass[aop,MSNbibl,nameyear,noautosecdot,dvips]{arximspdf}
\usepackage{graphicx}

\doi{10.1214/11-AOP723} %kopijuoti is PTS
\volume{41}
\issue{2}
\pubyear{2013}
\firstpage{744}
\lastpage{784}

\makeatletter
\newcommand{\ball}{\mathcal{B}}
\newcommand{\cyl}{\mathcal{C}}
\def\L{\mathcal{L}}
\def\ol{\overline}
\def\prt{\partial}

\renewcommand{\d}{\mathrm{d}}
\newcommand{\dist}{\operatorname{dist}}
\newcommand{\diam}{\operatorname{diam}}
\newcommand{\grad}{\operatorname{grad}}
\newcommand{\intrinsicdist}{\operatorname{dist}_{\mathrm{intr}}}
\newcommand{\length}{\operatorname{len}}
\newcommand{\origin}{\mathbf{o}}
\newcommand{\CAT}[1]{\operatorname{CAT}({#1})}

\newcommand{\Prob}[1]{\mathbb{P}[#1]}
\newcommand{\Reals}{\mathbb{R}}

\def\prt{\partial}
\def\wt{\widetilde}
\def\NN{\mathcal{N}}
\def\FF{\mathcal{F}}
\def\eps{\varepsilon}
\def\JJ{\mathbb{J}}
\def\KK{\mathbb{K}}

\newtheorem{thmm}{Theorem}
\newtheorem{cor}[thmm]{Corollary}
\newtheorem{prop}[thmm]{Proposition}
\newtheorem{lem}[thmm]{Lemma}
\newproclaim{defn}[thmm]{Definition}
\newproclaim{rem}{Remark}%[section]
\newtheorem{conjecture}{Conjecture}

\def\idist{\intrinsicdist}
\def\bX{\mathbf{W}}
\def\bY{\mathbf{Z}}
\def\TT{\mathcal{T}}
\newcommand{\eqref}[1]{(\ref{#1})}
\makeatother

\begin{document}
\begin{frontmatter}

\title{Shy couplings, $\operatorname{CAT}(0)$ spaces, and the lion and man\thanksref{T1}}
\thankstext{T1}{Supported in part by NSF Grants CCF-07-29537 and DMS-09-06743, and
by Grant N N201 397137, MNiSW, Poland.}
\runtitle{Shy couplings, $\CAT0$, and lion and man}

\begin{aug}
\author[A]{\fnms{Maury} \snm{Bramson}\ead[label=e1]{bramson@math.umn.edu}},
\author[B]{\fnms{Krzysztof} \snm{Burdzy}\ead[label=e2]{burdzy@math.washington.edu}}
\and
\author[C]{\fnms{Wilfrid} \snm{Kendall}\corref{}\ead[label=e3]{w.s.kendall@warwick.ac.uk}}
\runauthor{M. Bramson, K. Burdzy and W. Kendall}
\affiliation{University of Minnesota, University of Washington and
University~of~Warwick}
\address[A]{M. Bramson\\
School of Mathematics\\
University of Minnesota\\
Vincent Hall, 206 Church St. SE.\\
Minneapolis, Minnesota 55455\\
USA\\
\printead{e1}} %adresu isvedimo komanda gale!
\address[B]{K. Burdzy\\
Department of Mathematics\\
University of Washington\\
Box 354350\\
Seattle, Washington 98195\\
USA\\
\printead{e2}}
\address[C]{W. Kendall\\
Department of Statistics\\
University of Warwick\\
Coventry CV4 7AL\\
United Kingdom\\
\printead{e3}}
\end{aug}

% HISTORY:
\received{\smonth{7} \syear{2010}}
\revised{\smonth{7} \syear{2011}}

% ABSTRACT
%
\begin{abstract}
Two random processes $X$ and $Y$ on a metric space are said to be
\textit{$\varepsilon$-shy coupled} if there is positive probability of them
staying at least a positive distance $\varepsilon$
apart from each other forever.
Interest in the literature centres on nonexistence results subject to
topological and geometric
conditions; motivation arises
from the desire to gain a better understanding of probabilistic coupling.
Previous nonexistence results for co-adapted shy coupling of reflected Brownian
motion required convexity conditions; we remove these conditions by
showing the nonexistence
of shy co-adapted couplings of reflecting Brownian motion in any bounded
$\CAT0$ domain with boundary
satisfying uniform exterior sphere and interior cone conditions, for example,
simply-connected bounded planar domains with $C^2$ boundary.

The proof uses a
Cameron--Martin--Girsanov argument, together with a continuity property
of the
Skorokhod transformation and properties of the intrinsic metric of the domain.
To this end, a generalization of Gauss' lemma is established that shows
differentiability of the intrinsic distance
function for closures of $\CAT0$ domains with boundaries
satisfying uniform exterior sphere and interior cone conditions.
By this means, the shy coupling question is converted into a Lion and Man
pursuit--evasion problem.
\end{abstract}

% KEYWORDS
%
\begin{keyword}[class=AMS]
\kwd{60J65}.
\end{keyword}
\begin{keyword}
\kwd{CAT(0)}
\kwd{CAT($\kappa$)}
\kwd{co-adapted coupling}
\kwd{coupling}
\kwd{eikonal equation}
\kwd{first geodesic variation}
\kwd{Gauss' lemma}
\kwd{greedy strategy}
\kwd{intrinsic metric}
\kwd{Lion and Man problem}
\kwd{Lipschitz domain}
\kwd{Markovian coupling}
\kwd{pursuit--evasion problem}
\kwd{reflected Brownian motion}
\kwd{Reshetnyak majorization}
\kwd{shy coupling}
\kwd{Skorokhod transformation}
\kwd{total curvature}
\kwd{uniform exterior sphere condition}
\kwd{uniform interior cone condition}.
\end{keyword}\vspace*{-3pt}

\end{frontmatter}
%

%s1 ###
\section{\texorpdfstring{Introduction.}{Introduction}}\label{secintro}\vspace*{-3pt}
%s1.1 ###
\subsection{\texorpdfstring{Results and motivation.}{Results and motivation}}\label{secresults-and-motivation}

Benjamini, Burdzy and Chen (\citeyear{BenjaminiBurdzyChen-2007}) introduced the notion of \textit{shy
coupling}: a coupling of Brownian motions $X$ and $Y$ (more generally, of two\vadjust{\goodbreak}
random processes $X$ and $Y$ on a metric space) is said to be \textit{shy} if there is an
$\varepsilon>0$ such that $\Prob{\dist(X(t),Y(t))\geq\varepsilon\mbox{ for
all }t}>0$. For example consider Brownian motion $X$ on the circle: if
$Y$ is produced from $X$ by a nontrivial rotation then $X$ and $Y$ exhibit a shy
coupling, since $\dist(X,Y)$ is then constant. Interest in the
existence or nonexistence of such couplings arises from the study of couplings of reflected
Brownian motions, which occur in various contexts.
\citet{BenjaminiBurdzyChen-2007} discussed
existence and nonexistence of shy couplings for Brownian motions on
graphs and
for reflected Brownian motions in domains (connected open subsets of
Euclidean space) satisfying suitable boundary regularity conditions.
They restricted attention to Markovian couplings and we will
do essentially the same, by restricting attention to co-adapted
couplings. (This
is
only
slightly more general,
but
is more convenient for expression in terms of
stochastic calculus.) In particular
the results in
\citet{BenjaminiBurdzyChen-2007} showed
that no shy co-adapted
couplings can exist for reflected Brownian motion in convex bounded planar
domains with $C^2$ boundary satisfying a strict convexity condition (namely,
that the boundary contains no
% deletion of ``nontrivial'' suggested in final remarks by referee,
%received 4/10/2011
% nontrivial
line segments). Their argument used a large
deviations argument bearing some resemblance to
methods from
differential game
theory.
\citet{Kendall-2009a} showed that neither differentiability nor strict convexity
is required for the planar result, and also generalized the result to convex
bounded domains in higher dimensions whose boundaries need no longer be smooth
but still satisfy the regularity condition requiring triviality of all line
segments contained in the boundary. These more recent results are based
on direct proofs
using ideas from stochastic control.

The work described below both generalizes the above results and
also
shows that
absence of shyness is not confined to the case of convexity. We
consider a
bounded domain with boundary satisfying uniform exterior sphere and interior
cone conditions and that satisfies a $\CAT0$
condition (see Definition~\ref{m221})
when furnished with
the intrinsic metric, and we show that such domains cannot support shy
co-adapted couplings of reflected Brownian motions. We do this by
establishing a
rather direct connection between (the nonexistence of) Brownian shy co-adapted
couplings and deterministic pursuit--evasion problems. As part of this process,
we generalize Gauss' lemma (on the differentiability of the distance function)
to the case of closures of $\CAT0$ domains furnished with the intrinsic
metric and
satisfying uniform exterior sphere and interior
cone conditions.
It may not be evident to
the reader exactly how the stochastic and undirected notion of Brownian
motion can be
connected to the deterministic and intentional notion of a
pursuit--evasion problem, and it was not initially evident to us
[though, in
retrospect,
this
is latent in \citet{BenjaminiBurdzyChen-2007}], but
nonetheless the connection is both immediate and useful.

The pursuit--evasion problem in question is a well-known problem
concerning a Lion
chasing a Man in a disk, both travelling at unit speed:
R.~Rado's celebrated ``Lion and Man''
problem. Our shy coupling problem leads us to consider the generalization
in which the Lion chases the Man in a bounded domain which is $\CAT0$
in its
intrinsic metric. \citet{Isaacs-1965} is the classic reference for
pursuit--evasion problems;
\citet{Nahin-2007} provides an accessible exposition of the special
case of the
Lion and Man problem in the unit disk. Littlewood [(\citeyear{Littlewood-1986}), pages~114--117 in Bollobas'
extended edition] provides a brief description of the
Lion and Man problem with an indication of its history, including a
presentation of Besicovitch's celebrated proof that in the disc the Man can
evade the Lion indefinitely, even though the distance between Lion and
Man may
tend to zero. A generalization of discrete-time pursuit--evasion
to bounded $\CAT0$ domains is dealt with in
\citet{AlexanderBishopGhrist-2006}; we summarize
concepts from metric geometry and develop results required for the
continuous-time variant in Section~\ref{secdeterministic-pursuit}, and
it is here
that we generalize the Gauss lemma to the case of closures of $\CAT0$
domains with sufficient
boundary regularity (Proposition~\ref{propgauss}).

In particular, Section~\ref{secdeterministic-pursuit} rigorously
develops the geometric results required to reason with these
concepts in the context of the intrinsic metric for the domain~$D$
(determined by lengths of
paths restricted to lie within $D$). On a first reading one should feel
free to note only the general ideas of
Section~\ref{secdeterministic-pursuit}, and then to pass quickly on to
the probabilistic arguments
in the remaining sections of the paper.

In Section~\ref{seccat0-and-pursuit-evasion}, we describe how
continuous-time pursuit--evasion problems
can be solved in $\CAT0$ domains. We obtain an upper bound for the time
of $\varepsilon$-capture,
expressed in terms of domain geometry.
%KB new
Simultaneously with and independently of our research project,
Chanyoung Jun developed in his Ph.D. thesis [\citet{JunPhD}] a theory
of continuous pursuit in $\CAT\kappa$ spaces that overlaps somewhat
with our results.
%KB end new

Pursuit--evasion games involve control of the velocity of the pursuer
so as to
bring it arbitrarily close to the evader, regardless of what strategy
may be
adopted by the evader. In order to show nonexistence of Brownian shy couplings,
we investigate the possibility of bringing the Brownian pursuer (the
\textit{Brownian Lion}) arbitrarily close to the Brownian evader (the
\textit{Brownian Man}), regardless of how the Brownian motion of the
Brownian Man
is coupled to that of the Brownian Lion. The connection between
coupling and
deterministic Lion and Man problems is described in Section
\ref{seccoupling-to-pursuit}: a suitable pursuit strategy generates a
vector field $\chi$ on the configuration manifold generated by the
locations of
Brownian Lion and Man. (More pedantically, it generates a section of
the pullback of the tangent bundle of~$D$ to the configuration space of
the pursuer and
evader before capture.) If this pursuit strategy can be guaranteed to
bring the
Lion within $\varepsilon/2$ of Man by a bounded time $t_c$ in the
deterministic problem, then a Cameron--Martin--Girsanov argument
together with a
continuity property for the Skorokhod transformation shows that the Brownian
Lion has a positive probability of getting within distance $\varepsilon
$ of the Brownian Man,
whatever coupling strategy might be adopted by the Brownian Man.

The paper concludes with Section~\ref{seccomplements}, which discusses
possible extensions of these results, further questions, and conjectures.

%f1 ###
\begin{figure}

\includegraphics{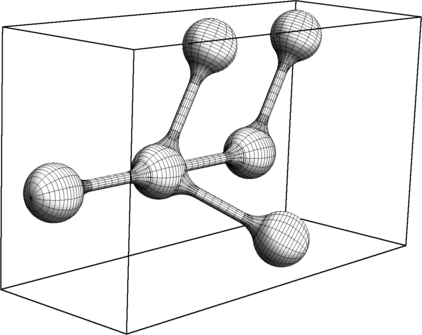}

\caption{A $\CAT0$ example which is the union of five dumbbells.}
\label{figdumbbells}
\end{figure}

We now state
the main results of this paper, using terms defined in Section~\ref{secdeterministic-pursuit}. Here and elsewhere in the paper,
we consider only domains in Euclidean space of dimensions $2$ or higher.
\begin{thmm}\label{thmno-shy-coupling}
Suppose that $D$ is a bounded domain with boundary satisfying uniform exterior
sphere and interior cone conditions, and which is $\CAT0$ in its intrinsic
metric. There can be no shy co-adapted coupling for reflected Brownian
motion in
$D$.
\end{thmm}

Examples of $\CAT0$ domains include convex domains and domains that are
the unions of a pair
of convex domains. See, for instance, \citet{BridsonHaefliger-1999} and
\citet{AlexanderBishopGhrist-2006},
where more general examples are also provided; in particular, a large
range of examples follows from
iterated application of
the result that if two $\CAT0$ domains have a geodesically convex
intersection then their union is $\CAT0$.
The exterior sphere and interior cone conditions in the theorem are
required in order to apply the results
of \citet{Saisho-1987} to generate reflected
diffusions using the Skorokhod transformation.

% Start of WSK insertion of Burdzy material, 08.10.10
The three-dimensional domain in Figure~\ref{figdumbbells} is $\CAT0$.
There are two different ways to see this.
First, it is easy to see that for every point on the boundary of the
domain, at most one of the principal curvatures is negative.
An alternative way to see that the domain is $\CAT0$ is to observe that
a single dumbbell (the union of two spheres and the connecting tube) is
a $\CAT0$ domain.
The whole set is the union of five dumbbells.
The nonempty intersections of the dumbbells are balls.

Remarkably, all bounded simply-connected planar domains are $\CAT0$ in their
intrinsic metrics. Thus, in the planar case, there is an immediate
consequence of Theorem~\ref{thmno-shy-coupling} which is
a strikingly powerful result
depending
principally
on topological conditions:
\begin{thmm}\label{thmno-shy-coupling2}
Suppose that $D$ is a simply-connected bounded planar domain with boundary
satisfying uniform exterior sphere and interior cone conditions. There
can be no
shy co-adapted coupling for reflected Brownian motion in $D$.
\end{thmm}

%s1.2 ###
\subsection{\texorpdfstring{Some basic tools for probabilistic coupling.}{Some basic tools for probabilistic coupling}}\label{secbasic-tools}
All probabilistic couplings considered here are \textit{co-adapted couplings},
which are defined
for general Markov processes in \citet{Kendall-2009a}. In essence, a
co-adapted coupling
of two Markov processes is a construction of the two Markov processes
on the same
probability space, which are adapted to the same filtration such that
each process possesses
the prescribed transition functions \textit{with respect to the common
filtration}.

In this paper, it suffices to work
with co-adapted couplings of $d$-dimensional Brownian motions: $B$ and
$\widetilde{B}$ are
said to be \textit{co-adaptively coupled Brownian motions} if they are
defined on the same probability space and adapted to the same filtration
$\{\mathcal{F}_t\dvtx t\geq0\}$ and if, in addition, both satisfy an independent
increments property taken \textit{with respect to the common filtration}:
\begin{eqnarray*}
&B_{t+s}-B_t \mbox{ is independent of }\mathcal{F}_t\mbox{ for all }t,
s\geq0 ,&
\\
&\widetilde{B}_{t+s}-\widetilde{B}_t \mbox{ is independent of }\mathcal
{F}_t\mbox{ for all }t, s\geq0 .&
\end{eqnarray*}
Note that $B_{t+s}-B_t$ and $\widetilde{B}_{t+s}-\widetilde{B}_t$ need
not be independent of each other.
\citeauthor{Kendall-2009a} [(\citeyear{Kendall-2009a}), Lemma 6] shows that one may represent such a
coupling using stochastic calculus,
possibly at the cost of augmenting the filtration by adding a further
independent Brownian motion $C$: there exist
$(d\times d)$-matrix-valued predictable random processes $\mathbb{J}$
and $\mathbb{K}$ such that
\[
\wt B
= \int\mathbb{J}^\top\,\d B + \int\mathbb{K}^\top\,\d C ;
\]
moreover, one may choose $\mathbb{J}^\top\mathbb{J}+\mathbb{K}^\top
\mathbb{K}$ to be equal to the $(d\times d)$
identity matrix at all times.

A pair of processes $X$ and $\widetilde{X}$ is said to form a
co-adapted coupling
if they can be defined by strong solutions of stochastic differential
equations driven
by $B$, $\widetilde{B}$, respectively. In the paper, we will employ the
stochastic differential
equation obtained from the Skorokhod transformation for reflected
Brownian motion in a domain $D$
of suitable boundary regularity, such as under uniform exterior sphere
and uniform interior cone conditions,
as discussed in Section~\ref{secdeterministic-pursuit}.
For $r>0$, set $\NN_{x,r} = \{\nu\in\Reals^d\dvtx  |\nu|=1, \ball(x+ r \nu
, r) \cap D = \varnothing\}$. The vectors $\nu$ can be
be viewed as ``exterior normal unit vectors at $x\in\partial D$''; note
that there may be more than one such vector
at a particular point $x\in\partial D$. The set $\NN_{x,r}$ is
decreasing in $r$, and the uniform exterior sphere condition
asserts that $r$ can be chosen so that, for all $x\in\prt D$,
$\NN_{x,r} \ne\varnothing$, with $\NN_{x,r}=\NN_{x,s}$ for $0<s\leq r$.
Under uniform exterior sphere and uniform interior cone conditions,
\citet{Saisho-1987} has shown that, given
a driving Brownian motion $B$,
there exists a unique solution pair $(X, L^X)$ satisfying
%
%e1 ###
\begin{eqnarray*}
&\d X = \d B - \nu_X \,\d L^X ,&
\\
&L^X\mbox{ is nondecreasing and increases only when }X\in\partial D
,&\\
&\nu_X \in\NN_{X,r} .&
\end{eqnarray*}
Thus $L^X$ may be viewed as the local time of the reflected Brownian
motion $X$ on the boundary $\partial D$.

In this paper, all vectors
are assumed to be column vectors unless specified otherwise.

% \section[CAT(0) geometry and the deterministic pursuit-evasion
% problem]{$\CAT0$ geometry and the deterministic pursuit-evasion
% problem}\label{secdeterministic-pursuit}
%s2 ###
\section{\texorpdfstring{$\CAT0$ geometry and the deterministic pursuit--evasion
problem.}{$\CAT0$ geometry and the deterministic pursuit--evasion
problem}}\label{secdeterministic-pursuit}

Recall that the \textit{intrinsic metric} for a domain $D$ is generated
by the
infimum of Euclidean lengths $\length(\gamma)$ of smooth
connecting paths $\gamma$ lying wholly within the domain.
(The definition is typically formulated in the context of
general metric spaces and regularizable paths.)

\begin{defn}
The \textit{intrinsic distance} between two points $x$ and $y$ in a domain~$D$ is given by
%
%e2 ###
\begin{equation}\label{eqintrinsic-metric}
\intrinsicdist(x,y) =
\inf\{\length(\gamma)\dvtx \gamma\mbox{ is a smooth path connecting }x\mbox{
and }y\mbox{ in }D\} .\hspace*{-30pt}
\end{equation}
\end{defn}

For a domain $D$, a standard compactness argument shows that paths
attaining the
infimum of \eqref{eqintrinsic-metric} will always exist in the closure
of the
domain: these are called \textit{intrinsic geodesics}.

% WSK July 2011: unify CAT definitions, using
As described in
\citeauthor{BridsonHaefliger-1999} [(\citeyear{BridsonHaefliger-1999}), Section~II.1, Definition
1.1] [see also \citet{BuragoBuragoIvanov-2001}], one can define simple curvature
conditions for metric spaces such as $(D, \intrinsicdist)$, based on the
behaviour of geodesic triangles. We first give the case of comparison
with flat Euclidean space
(which has zero curvature).
\begin{defn}\label{m221}
We say that $(D, \intrinsicdist)$ is
% adding ''a ... domain'' suggested in final remarks by referee,
%received 4/10/2011
a
$\CAT0$
domain
% ends
if the following triangle
comparison holds:
Suppose that $\Gamma_{a,b}$, $\Gamma_{a,c}$ and $\Gamma_{b,c}$ are
unit-speed intrinsic geodesics for $D$, connecting points $a$ to $b$,
$a$ to $c$ and $b$ to $c$, respectively. Then, for all such geodesic
triangles,
\[
\intrinsicdist(\Gamma_{a,b}(s),\Gamma_{a,c}(t)) \leq r(s,t) ,
\]
where $r(s,t)$ is the distance between points at distance $s$,
respectively, $t$, from
$\widetilde{a}$ along the side $\widetilde{a}\widetilde{b}$,
respectively,
$\widetilde{a}\widetilde{c}$, of an ordinary Euclidean triangle
$\widetilde{a}\widetilde{b}\widetilde{c}$ that has the same side lengths.
\end{defn}
%
% \marginpar{WSK: Label edges in original triangle.}
%
%f2 ###
\begin{figure}

\includegraphics{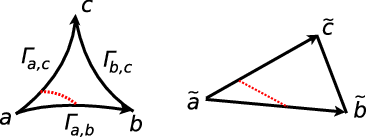}

\caption{Illustration of the $\CAT0$ condition.}
\label{figshy2-cat0}
\end{figure}

Consequently, chords of triangles in $(D, \intrinsicdist)$ are shorter
than comparable chords
of the comparable Euclidean triangles, as illustrated in Figure~\ref{figshy2-cat0}.

%KB new
\citeauthor{BridsonHaefliger-1999} [(\citeyear{BridsonHaefliger-1999}), Section~II.1, Definition 1.1] actually
introduces the more general notion
of a $\CAT\kappa$
domain [see also \citet{AlexanderBishopGhrist-2009}, Appendix~A]. Here we describe the case when
comparisons are drawn with triangles on a sphere of radius $1/\sqrt
{\kappa}$, for $\kappa>0$
(hence the sphere has curvature $\kappa$). It is necessary here to
restrict attention to suitably small triangles, as measured by perimeter.
\begin{defn}\label{july11}
We say that $D$ is a $\CAT\kappa$
domain for $\kappa> 0$ if any two distinct points with distance less
than $\pi/\sqrt{\kappa}$ are joined by a geodesic and the distance
between any two
points of any geodesic triangle $\triangle pqr$ of perimeter less than
$2\pi/\sqrt{\kappa}$
is no greater than the distance
between the corresponding points of the model triangle $\triangle\wt p
\wt q \wt r$ with the same sidelengths in
the 2-dimensional Euclidean sphere of radius $1/\sqrt{\kappa}$.
\end{defn}
%
%KB end new

\begin{rem*}
% $\text{CAT}$ is an acronym introduced by Gromov,
% . It stands for
Gromov introduced the acronym $\operatorname{CAT}$, standing for
{C}artan, {A}leksandrov, {T}oponogov.
% In the notation $\CAT0$, the zero stands
% for the zero curvature of the Euclidean space within which the
%comparison triangle is situated.
%KB start new
In this paper, we will mostly consider spaces $\CAT\kappa$ with $\kappa=0$.
% We will not suppress the redundant parameter 0 in $\CAT0$ to make it
%easier for the reader to compare our results in the present paper with
%those in
We include some results concerning the $\CAT\kappa$ case with $\kappa
>0$ because they will be used in
the forthcoming paper \citet{BBK}.
% , where $\CAT\kappa$ domains will be used.
%KB end new
% WSK July 2011: end unification edits.
\end{rem*}

\begin{rem*}
% WSK July 2011: angles for CAT(k) edits.
As noted in \citeauthor{BridsonHaefliger-1999} [(\citeyear{BridsonHaefliger-1999}), Proposition~II.3.1] [see also
\citet{BuragoBuragoIvanov-2001}, Section~4.3], in $\CAT\kappa$ spaces
% (more
% generally, Alexandrov spaces of nonpositive curvatures),
the notion of angle is
well-defined for (locally) minimal geodesics.
% WSK July 2011: end angles for CAT(k) edits.
\end{rem*}

Consequently, geodesics in a $\CAT0$ space diverge at %KB least
most
as fast as corresponding
geodesics in Euclidean space. Note that $\CAT0$ is a global condition,
applying to all possible geodesic triangles. In particular it can be
shown that
$\CAT0$ spaces are always simply-connected and indeed contractible
[\citet{BridsonHaefliger-1999}, Proposition~II.1.4, or \citet{AlexanderBishopGhrist-2009}, Appendix~A].\vadjust{\goodbreak}

\begin{figure}

\includegraphics{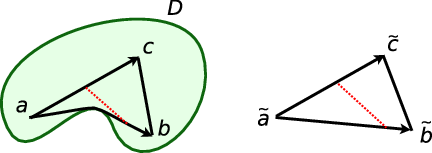}

\caption{Illustration of the $\CAT0$ property for a bounded simply-connected
planar domain. The effect of the boundary is to make the triangle ``skinnier''
than its Euclidean counterpart, thus establishing the $\CAT0$ comparison
property.}
\label{figshy2-cat0-D}\vspace*{-3pt}
\end{figure}

Remarkably, bounded \textit{planar} domains are $\CAT0$ if they are
simply-connected; see \citet{Bishop-2008} for a careful proof. Readers may
convince themselves of this at an intuitive level by drawing pictures
(as exemplified in Figure~\ref{figshy2-cat0-D}); as is the case with other
foundational results in metric spaces, the rigorous proof requires delicate
reasoning.
%f3 ###

We now introduce two complementary notions of boundary regularity following
\citet{Saisho-1987}. An \textit{exterior sphere condition} (also called
\textit{weak
convexity}) requires that every boundary point is touched by at least one
external sphere.
Here and in the following, let $\ball(y,s)$ denote the open Euclidean
ball of radius $s$ centered on $y$.\vspace*{-2pt}
\begin{defn}[{[Uniform exterior sphere condition, from \citet{Saisho-1987}, Section~1, Condition $(A)$]}]\label{defexterior-sphere-condition}
%KB next line is for typesetting reasons
A domain $D$ is said to satisfy a \textit{uniform exterior sphere condition,
based on radius $r$} if, for every $x\in\partial D$, the set of
``exterior normals''
$\NN_{x,r} = \{\nu\in\Reals^d\dvtx  |\nu|=1, \ball(x+ r \nu, r) \cap D =
\varnothing\}$
is nonempty, with $\NN_{x,r}=\NN_{x,s}$ for $0<s\leq r$.\vspace*{-2pt}
\end{defn}

Thus a uniform exterior sphere condition allows one to move a
fixed ball all the way around the outside of the domain
boundary. In particular, $D$ can have no ``inward-pointing
corners''. Here is a simple observation which will be useful
later and corresponds to the intuition about being able to move
a fixed ball about $D$; such $D$ may be represented
as intersections of complements of balls,
in a manner entirely analogous to the
representation of a convex set as the intersection of
half-planes (so justifying the alternative term ``weak
convexity'').
\begin{lem}\label{lemsg-for-uesc}
Suppose that the domain $D$ satisfies a uniform exterior sphere condition
based on radius $r$. Then
\[
\overline{D} = \bigcap\{\ball(z,r)^c \dvtx  \ball(z,r)\cap
D=\varnothing\} .
\]
\end{lem}

\begin{pf}
Let the Minkowski sum $A\oplus B$ of two Euclidean sets $A$ and $B$
be $A\oplus B=\{x+y\dvtx x\in A, y\in B\}$.
Certainly $F=\bigcap\{\ball(z,r)^c \dvtx  \ball(z,r)\cap
D=\varnothing\}$ is closed, since $\ball(z,r)$ is an open ball.
Moreover $D\subseteq F$; hence $\overline{D}\subseteq F$.
Furthermore $F\subseteq\overline{D}\oplus\ball(\origin,r)$, where
$\origin$ is the origin of the ambient Euclidean space.\looseness=-1\vadjust{\goodbreak}

Following \citeauthor{Saisho-1987} [(\citeyear{Saisho-1987}), Remark 1.3], because of the uniform exterior
sphere condition, we can define a projection $x\mapsto\overline{x}$ from
$\overline{D}\oplus\ball(\origin, r)$ onto $\overline{D}$ using the
Euclidean metric.
Consider $x\in\overline{D}\oplus\ball(\origin,r)$. Then the projection
$\overline{x}\in\overline D$
is defined; moreover,
if $\overline{x}\in\partial D$ and $x \notin\overline D$, then
\[
\frac{\overline{x}-x}{|\overline{x}-x|} \in
\mathcal{N}_{\overline{x},r}
\]
is a unit vector whose offset produces a tangent sphere of radius $r$ at
$\overline{x}$ [using the argument of \citet{Saisho-1987}]. But this implies
that
if $x\in(\overline{D}\oplus\ball(\origin,r))\setminus D$ then
\[
x \in
\ball\biggl(\overline{x}-r\frac{\overline{x}-x}{|\overline{x}-x|},r\biggr)
\]
and so $x\notin F$. Accordingly $\overline{D}=F$ as required.
\end{pf}

On the other hand, a \textit{uniform interior cone condition} requires
that any boundary point
supports a bounded cone truncated to the boundary of a ball, and
moreover that
the cone may be translated locally within the
domain.

\begin{defn}[{[Uniform interior cone condition, from \citet{Saisho-1987}, Section~1,
Condition $(B^\prime)$]}]\label{defUICC}
A domain
$D$ is said to satisfy a \textit{uniform interior cone
condition, based on radius $\delta>0$ and angle
$\alpha\in(0,\pi/2]$,} if, for every $x\in
\partial D$, there is at least one unit vector $\mathbf{m}$ such that
the cone
$C(\mathbf{m})=\{z\dvtx \langle z,\mathbf{m}\rangle> |z|\cos\alpha\}$ satisfies
\[
\bigl(y + C(\mathbf{m})\bigr)\cap\ball(x,\delta) \subseteq D
\qquad\mbox{for all } y\in D\cap\ball(x,\delta) .
\]
We say that
% the cone $y+C(m)$
the cone $y+C(\mathbf{m})$
\textit{is
based on $y$
and angle $\alpha\in(0,\pi/2]$}.
\end{defn}

Thus a uniform interior cone condition implies that the ``outward-pointing
corners'' must not be too sharp. Note that \citeauthor{Saisho-1987} actually
uses a slightly weaker condition with less intuitive content [\citet{Saisho-1987}, Condition
$(B)$]; we do not consider this weaker notion further in what
follows.

In fact, the property of a domain satisfying a uniform interior cone
condition is equivalent to it being a Lipschitz domain.
\begin{defn}[(Lipschitz domain)] \label{deflipschitz-domain}
Recall
that a function $f\dvtx  \Reals^{d-1} \to\Reals$ is
\textit{Lipschitz, with constant $\lambda< \infty$}, if $|f(x) - f(y)|
\leq\lambda
|x-y|$ for all
$x,y \in\Reals^{d-1}$.
A domain $D$ is said to be \textit{Lipschitz, with
constant $\lambda$,} if there exists $\delta>0$ such that, for
every $x \in\prt D$, there exists an orthonormal
basis $e_1, e_2, \ldots, e_d$
and a Lipschitz function $f\dvtx  \Reals^{d-1} \to
\Reals$, with constant $\lambda$, such that
\[
\ball(x,\delta) \cap D =
\{y\in\ball(x,\delta)\dvtx  f(y_1, \ldots, y_{d-1})<y_d\} ,
\]
where we write $y_1=\langle y,e_1\rangle, \ldots, y_d=\langle
y,e_d\rangle$.
\end{defn}

The equivalence of Definitions~\ref{defUICC} and \ref
{deflipschitz-domain} depends on the fact that the cone axis vector in
Definition~\ref{defUICC} is chosen to be the same for all $y\in\ball
(x,\delta)$,
and so can be used as $e_d$ in the orthonormal basis for $\ball(x,\delta
)$ required in Definition~\ref{deflipschitz-domain}.
The constants $\lambda$ and $\alpha$ in Definitions~\ref{defUICC}
and~\ref{deflipschitz-domain} are related by $\lambda= \cot\alpha$,
while the two~$\delta$'s of Definitions~\ref{defUICC} and \ref
{deflipschitz-domain} may be taken to be equal.
Note too that if the uniform interior cone/Lipschitz domain property
holds for a given $\delta>0$, then evidently it also holds for all
smaller $\delta$.

If a domain satisfies a uniform interior cone condition, then the
intrinsic metric and Euclidean metric properties are closely related.
\begin{lem}\label{lemdiameter}
A domain $D$ that is bounded in Euclidean metric and
satisfies a uniform interior cone condition must have finite
intrinsic diameter.
\end{lem}

% WSK correction and clarification of proof 5 July 2011
%
\begin{pf}
Certainly $\intrinsicdist(x,y)$ is a
continuous
function of $(x,y)$ in the open set $ D\times D$
and takes only finite values there.
% The
Note that the
domain $D$ is path-connected, being an
open connected subset
of Euclidean space.

% \marginpar{WSK: $\frac12\delta\sin\alpha$ I think! Also, $\delta/2$
%not $\delta$.}
Suppose that $D$ satisfies a uniform interior cone condition based on radius
$\delta>0$ and angle $\alpha\in(0,\pi/2]$.
If $\textbf{m}$ is a unit vector for the interior cone condition at
$x\in\partial D$
then geometrical arguments show that $x+\frac12\delta\textbf{m}$ is at
least $\frac12\delta\sin\alpha$
from the exterior $D^c$.
Choosing
$\delta^\prime$ with
% $0<\delta^\prime<\delta\tan\alpha$,
$0<\delta^\prime<\frac12\delta\sin\alpha$, it follows that any such
$x+\frac12\delta\textbf{m}$
belongs to
\[
D\ominus\ball(\origin,\delta^\prime) \stackrel{\mathrm{def}}{=}
\bigl(D^c\oplus\ball(\origin,\delta^\prime)\bigr)^c
=
\{x\in D\dvtx
{\ball(x,\delta^\prime)}\subset D\} ,
\]
which itself % is
is closed. Inheriting boundedness from $D$, it is therefore compact in the
Euclidean topology, and hence also in the topology derived from the intrinsic
metric, since the two metrics are locally equal away from the boundary
of~$D$.
Hence $\{\intrinsicdist(x,y)\dvtx  x,y\in D\ominus\ball(\origin,\delta^\prime
)\}$
attains a maximum value, which is therefore finite. However, for \textit{any}
$x^\prime,y^\prime\in D$, we have
%
%e3 ###
\begin{equation}\label{eqdiameter}
\intrinsicdist(x^\prime,y^\prime) \leq %2
\delta+
\sup\{\intrinsicdist(x,y)\dvtx  x,y\in D\ominus\ball(\origin,\delta^\prime
)\} ,
\end{equation}
because we have used
the uniform interior cone condition
% assures us
to ensure % that
that from each
point on the boundary there is a straight-line segment of
length $\frac12\delta$ to $\{\intrinsicdist(x,y)\dvtx  x,y\in
D\ominus\ball(\origin,\delta^\prime)\}$. Hence the
intrinsic diameter
must be bounded by the right-hand side of \eqref{eqdiameter}.
\end{pf}
%
% WSK end of correction and clarification of proof 5 July 2011

The full force of the uniform interior cone condition is not required
for the
above result;
the proof does not require coordination of the
directions of interior cones at different base-points.
The full force of
the uniform interior cone condition
assures us that any path
of finite length leading in $D$ to a point $x$ on the boundary of $D$ can
be deformed continuously in $D$ into one which in its final phase is
the~segment on which the interior cone at $x$ is based. Moreover, the lengths
of the curves throughout this deformation can be constrained to be arbitrarily
close to the length of the original path. This allows us to view $D$ as a
topological manifold with boundary, which is continuously embedded in
the ambient
Euclidean space. More than this, it shows that the completion $\hat{D}$ of
$D$ under the intrinsic metric can be identified with the Euclidean closure
$\overline{D}$ and moreover that the intrinsic metric and the Euclidean metric
actually endow $\overline{D}$ with the same topology. Finally,
\citeauthor{BridsonHaefliger-1999} [(\citeyear{BridsonHaefliger-1999}), Corollary II.3.11] show that the closure
$\overline{D}$, viewed as the completion $\hat{D}$
of $D$ in intrinsic metric, inherits $\CAT0$ structure from $D$.

%s2.1 ###
\subsection{\texorpdfstring{Regularity for geodesics.}{Regularity for geodesics}}
We
wish to
consider pursuit--evasion in a bounded $\CAT0$ domain. Lion and Man both
move with unit speed, with the Lion seeking to draw closer to the Man
by using a
``greedy'' pursuit strategy (which is not necessarily optimal).
This Lion strategy can be phrased in terms of an $\Reals^d$-valued
field $\chi$ of unit vectors defined on the configuration space
$(\overline{D}\times
\overline{D})\setminus\{(x,x)\dvtx x\in\overline{D}\}$, such that $\chi(x,
y)$ is the initial velocity of the
unit-speed geodesic moving from $x$ to $y$. (This is the vector field
described pedantically in Section~\ref{secintro} as a section of the
pullback of the tangent bundle of $D$ to the configuration space of the
pursuer and
evader before capture.)

% WSK: 5 July 2011 Needed to change verbal summary, which wasn't quite
%right!
We first show that the combination of uniform exterior sphere and
uniform interior cone/Lipschitz conditions
implies that, working locally,
every boundary point of the intersection of the domain $D$ with a
suitable $2$-plane will support an
exterior sphere, albeit with smaller radius.
% We first show that the combination of uniform exterior sphere and
%uniform interior cone / Lipschitz conditions
% implies that every boundary point of the intersection of the domain
%$D$ with a suitable $2$-plane will support an
% exterior sphere, albeit with smaller radius.
%
\begin{lem}\label{lemplanar-UESC}
Suppose that $D$ is a domain satisfying a uniform exterior sphere
condition based on radius $r>0$, and a uniform interior cone
condition based on radius $\delta>0$ and angle
$\alpha\in(0,\pi/2]$.
Suppose that $z\in\partial D$ and $e_d$ is the $d$th vector in
the orthonormal basis
corresponding to $z$ as in Definition~\ref{deflipschitz-domain}.
Let $P$ be a $2$-plane intersecting $D$ and containing $z$ and $z+e_d$.
Then there exists $w\in P$, with $|w-z|=\dist(w,D\cap P)=r \sin\alpha$ and
%$\ball(w,s)\cap(D\cap P)=\varnothing$.
$\ball(w,r \sin\alpha)\cap(D\cap P)=\varnothing$.
\end{lem}

% WSK3: the modification here is required because the true exterior
%sphere condition requires a stabilization
% property for the ensemble of exterior spheres at each point, which is
%not required for $D\cap P$. (Nor do
% I see a way to prove it, if indeed it is true at all!)
% Note that $z\in\partial(D\cap P)$ implies $z\in\partial D$.
Since the interior cone condition is uniform, the lemma shows that
every point in the boundary of $D\cap P$
near $z$ must
support an exterior sphere of radius $r\sin\alpha$.
% WSK 5 July 2011 altering proof to localize it.
%
%
%f4 ###
\begin{figure}

\includegraphics{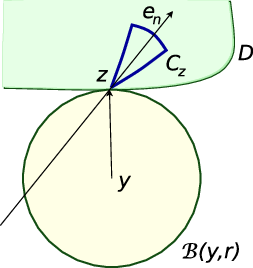}

\caption{Illustration of interior cone $C_z$ and exterior ball $\ball(y,r)$ at
$z\in D$.}
\label{figshy2-2d-cone}
\end{figure}

\begin{figure}[b]

\includegraphics{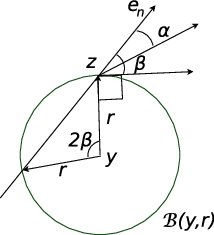}

\caption{Two-dimensional section of Figure \protect\ref{figshy2-2d-cone}
illustrating the underlying two-dimensional geometry.}
\label{figshy2-2d-cone-plan}
\end{figure}

\begin{pf*}{Proof of Lemma \protect\ref{lemplanar-UESC}}
Suppose that $z\in\partial D$ with $e_d$ defined as above.
Let~$P$ be a $2$-plane containing $z$ and $z+e_d$.
% Fix $z\in\partial D$ and choose a $2$-plane $P$ containing $z$
Since $z\in\partial D$, there is an exterior sphere touching $z$,
defined by a ball $\ball(y,r)\subseteq D^c$ with $z\in\overline{\ball(y,r)}$.
By
Definition~\ref{deflipschitz-domain}, the cone
\[
C_z = \{w\dvtx  \langle w-z,e_d\rangle> |w-z|\cos\alpha\}
\]
lies locally in $D$, in the sense that $C_z\cap\ball(z,\delta)\subseteq
D$ (see Figure~\ref{figshy2-2d-cone}).
If $\frac{\pi}{2}+\beta$ is the angle between\vadjust{\goodbreak} $e_d$ and $y-z$, then
two-dimensional geometry (Figure~\ref{figshy2-2d-cone-plan}) shows that
\[
\min\{|y-(\gamma e_d+z)|\dvtx \gamma\in\Reals\} = r\cos\beta.
\]
But $\beta\geq\alpha$ if $C_z\cap\ball(z,\delta)\subseteq D$ and $\ball
(y,r)\subseteq D^c$;
moreover, the line $\{
\gamma e_d+z\dvtx
\gamma\in\Reals\}$ must lie in $P$.
Hence the distance from $y$ to $P$ is at most
$r\cos\beta\leq r\cos\alpha$.
Consequently the radius of the disk $\ball(y,r)\cap P$ is at least $r
\sin\alpha$;
since $z\in\partial(D\cap P)$ and $\ball(y,r)\cap P$ is an exterior
sphere to $z$ in $P$,
the lemma follows.
\end{pf*}

%f5 ###

% WSK 5 July 2011 end of altering proof to localize it.

% WSK 5 July 2011: edited next sentence for style!
We can now establish some important technical
consequences of the uniform exterior sphere and interior cone conditions
% , employed
together
with
the
$\CAT0$
condition; namely, that
the Euclidean and intrinsic distances are locally comparable,
and that the vector field $\chi$ is continuous with reference to
the common topology of the Euclidean metric and the intrinsic metric,
and hence
is uniformly continuous over regions for which the two arguments are
well-separated. This is spelled out in the following proposition.
%KB new
% We state
% WSK 5 July 2011
In fact we state and prove a generalization of
the result for
$\CAT\kappa$ domains with $\kappa\geq0$, so that we can apply it in
the forthcoming paper \citet{BBK}.
%KB end new
%
\begin{prop}\label{proplipschitz}
Suppose that $D$ is
%KB a $\CAT0$ domain,
a $\CAT\kappa$ domain with $\kappa\geq0$,
bounded in the Euclidean metric
and satisfying a uniform exterior sphere
condition based on radius $r>0$, and a uniform interior cone
condition based on radius $\delta>0$ and angle
$\alpha\in(0,\pi/2]$.
%KB new
We can and will assume without loss of generality that $\lambda=\cot
\alpha> 1$.
%KB end new
%
\begin{longlist}[(1)]
\item[(1)]\hypertarget{itemlipschitzpart1-metric}
Suppose $a$, $b\in\overline{D}$ are close in the Euclidean metric, in
the sense that
%KB formula changed, not only reformatted
%KB $|a-b|<\min
%KB \{\delta,2r\sin\alpha\}$.
%
%e4 ###
\begin{equation}\label{july23}
|a-b|<\min\{\delta/(4\lambda),2r\sin\alpha\}.
\end{equation}
Then
%
%e5 ###
\begin{equation}
2r\sin\alpha \sin\biggl(\frac{\intrinsicdist(a,b)}{2r\sin\alpha}\biggr) \leq
|a-b| \leq \intrinsicdist(a,b) .
\label{eqlipschitzpart1-metric}
\end{equation}
\item[(2)]\hypertarget{itemlipschitzpart2-differentiable}
%KB new
Suppose that $\kappa=0$.
%KB end new
Intrinsic geodesics for $D$
[necessarily minimal, by the $\CAT0$ condition]
are
continuously differentiable and their direction fields satisfy a Lipschitz
property\vspace*{-1pt}
with constant $\frac{4}{\sqrt{3}}\frac{1}{2r\sin\alpha}$ that therefore
holds uniformly for all minimal intrinsic geodesics in $D$ and
hence in $\overline{D}$ [since $\CAT0$ geodesics depend continuously on
their endpoints].
%KB new
For $\kappa>0$, the same conclusion holds for minimal intrinsic
geodesics with endpoints in $\overline{D}$ which are separated by
intrinsic distance
strictly less than $\pi/\sqrt{\kappa}$.
%KB end new
%
\item[(3)]\hypertarget{itemlipschitzpart3-chi}
For $x$, $y$ in $\overline{D}$
%KB new
with $\intrinsicdist(x,y) < \pi/\sqrt{\kappa}$ (as usual, $\pi/\sqrt
{0}=\infty$),
%KB end new
let $\chi(x,y)$ be the unit vector at
$x$ pointing along the unique intrinsic geodesic $\gamma^{(x,y)}$ from
$x$ to $y$.
Then $\chi(x,y)$ depends continuously on
%KB $(x,y)\in
%KB (\overline{D}\times\overline{D})
%KB \setminus\{(x,x):x\in\overline{D}\}$
$(x,y)$ in
$ A=\{(x,y)\in\overline{D}\times\overline{D}\dvtx
0<\intrinsicdist(x,y) < \pi/\sqrt{\kappa}\}$
and hence is uniformly
continuous over
%KB regions for which the two arguments are well-separated.
compact subregions of $A$.
\end{longlist}
\end{prop}

\begin{pf*}{Proof of part (\hyperlink{itemlipschitzpart1-metric}{1})}
% WSK3: Exchanged $d$ and $\rho$ in this proof, to agree with
%convention in remainder of paper.
% I did this with careful search-and-replace: hopefully no errors have
%crept in?!
Definitions~\ref{defUICC} and~\ref{deflipschitz-domain} are equivalent,
so the domain $D$ is Lipschitz with constant
$\lambda=\cot\alpha$.
For each ball of radius $\delta$, we may therefore construct a
coordinate system $e_1, \ldots, e_d$ and a Lipschitz function $f$
to implement the Lipschitz property of $D$.

% WSK 7 July 2011 amended control on $|a-b|$.
Consider $a$, $b\in\overline{D}$ with $|a-b|<\min\{\delta/(4\lambda
),2r\sin\alpha\}$.
If the line segment $S$ between $a$ and $b$ does not intersect
$\partial D$,
then it must form the (unique, minimal) intrinsic geodesic between
$a$ and $b$, and \eqref{eqlipschitzpart1-metric} follows immediately.
If $S$ does not intersect $\operatorname{int}(D^c)$,
then we can cover the intersection $S\cap\partial D$ with
% WSK 7 July 2011 a single ball will suffice!
%finitely many balls
a single
$\ball(z,\delta)$
% , $z\in\partial D$,
(for $z\in\partial D$)
and use the unit vector $e_d$ corresponding to
% each
the
ball
(equivalently, the unit vector defining the cone for the ball)
to perturb $S$ to a regularizable path in $D$ (save for the endpoints)
with length arbitrarily close to that of $S$.
Hence $S$ is the intrinsic geodesic between
$a$ and $b$, and therefore \eqref{eqlipschitzpart1-metric} follows immediately.
So we can confine our attention to the case when $a\neq b$ and $S$
intersects $\operatorname{int}(D^c)$.

Applying Definition~\ref{deflipschitz-domain} to $\ball(a,\delta)$,
there is a Lipschitz function $f\dvtx \Reals^{d-1}\to\Reals$, with Lipschitz
constant $\lambda=\cot\alpha$,\vadjust{\goodbreak}
and an orthonormal basis $e_1, \ldots, e_d$, such that
%KB \[
%
%e6 ###
\begin{equation}\label{july21}
\ball(a,\delta)\cap D = \{y\in\ball(a,\delta)\dvtx f(y_1,\ldots
,y_{d-1})<y_d\},
\end{equation}
%
%KB \]
where $y_1=\langle y,e_1\rangle, \ldots, y_d=\langle y,e_d\rangle$.
%KB new
%KB \cyl is encoded as cal C (can be changed globally)
% Recall that $\lambda> 1$.
% WSK 7 July 2011 the region $\cyl(a)$ must be centred on $a$!
% Let
Consider
\begin{eqnarray*}%\label{july21}
\cyl(a) &=& \{y\in\Reals^d\dvtx  (|y_1-a_1|^2 + \cdots+
|y_{d-1}-a_{d-1}|^2)^{1/2} < \delta/(4\lambda) ,\\
&&\hspace*{183pt}\phantom{\{}  |y_d-a_d| < \delta
/2\}
\end{eqnarray*}
and note that, since $\lambda>1$, it is a consequence of
\eqref{july21} that $\cyl(a)\subset\ball(a,\delta)$.
% It follows from
% Moreover it is a consequence of
% \eqref{july21} that
% WSK 7 July 2011 we need to use the Lipschitz property to show
%boundary stays in $\cyl(a)$
Moreover~$f$ has Lipschitz constant $\lambda$, so we can control the
behaviour of that
part of the boundary of $D$ lying within $\cyl(a)$:
%
%e7 ###
\begin{eqnarray}\label{july22}
&&\{y\in\Reals^d\dvtx (|y_1-a_1|^2 + \cdots+ |y_{d-1}-a_{d-1}|^2)^{1/2} <
\delta/(4\lambda),\nonumber\\
&&\hspace*{150pt}\quad f(y_1,\ldots,y_{d-1})=y_d\}
\\
&&\qquad=
\{y\in\cyl(a)\dvtx f(y_1,\ldots,y_{d-1})=y_d\}
=
\prt D \cap\cyl(a) .\nonumber
\end{eqnarray}
%
%KB end new
%KB There are major (mostly unmarked) changes in the rest of the proof
%of part (1)
%KB needed to adjust the reasoning from the ball \ball(a,\delta)
%KB to the cylinder \cyl(a)
Applying Lemma~\ref{lemplanar-UESC} to the $2$-plane
\[
P = a + \mbox{linear span}\{b-a,e_d\} ,
\]
%
% Note that $x\in P$, and $D\cap P\neq\varnothing$.
% WSK 7 July 2011 need to take care here: we can't choose a single
%plane $P$ over all $D$, only within $\ball(a,\delta)$
every boundary point of $D\cap P\cap\cyl(a)$ supports an exterior disk
of radius $r\sin\alpha$.
[Note that the Lipschitz representation implies that $(\partial D)\cap
P\cap\cyl(a)=\partial(D\cap P)\cap\cyl(a)$.]
We shall use these exterior disks to construct a short path between $a$
and $b$.

It follows from \eqref{july23} and \eqref{july22} that the two rays
from $a$ and $b$ along the direction $e_d$ must lie in $P\cap D$ until
they leave $\cyl(a)$:
%
%e9 ###
%e8 ###
\begin{eqnarray}\label{eqrays}
\cyl(a) \cap\{a +\gamma e_d\dvtx  \gamma>0\}& \subseteq &P\cap D ,
\nonumber
\\[-8pt]
\\[-8pt]
\nonumber
\cyl(a) \cap\{b +\gamma e_d\dvtx  \gamma>0\}& \subseteq& P\cap D .
\end{eqnarray}
We set $u$, $v$ to be the intersections of these rays with $\partial\cyl(a)$.
For each
%KB \lambda is the Lipschitz constant; changed to \eta without marking
%KB $\lambda\in(0,1)$,
$\eta\in(0,1)$,
consider the point $\eta a + (1-\eta) b$ and the open segment which is
the intersection of the corresponding ray with $\cyl(a)$,
namely
\[
\cyl(a) \cap\{\eta a + (1-\eta) b +\gamma e_d\dvtx  \gamma>0\} .
\]
It follows from \eqref{july23} and \eqref{july22} that a nonempty
final sub-segment
\[
\cyl(a) \cap\{\eta a + (1-\eta) b +\gamma e_d\dvtx  \gamma>\gamma_\eta\}
\]
must lie in $D\cap P$. But then any exterior disk for $\eta a + (1-\eta
) b +\gamma_\eta e_d$ has to avoid the rays defined in \eqref{eqrays}
as well as the above nonempty final sub-segments; it must not
intersect the segments $[a,u]$ and $[b,v]$, and also may not intersect
that portion of $\partial\cyl(a)$ which intersects rays $\{\eta a +
(1-\eta) b +\gamma e_d\dvtx  \gamma>0\}$ (see Figure \ref
{figshy2-2d-plane}). Consequently (since $|a-b|<2r\sin\alpha$) such an
exterior disk
must have center lying on the side of the line through $a$ and $b$
which is opposite to the side containing $u$ and $v$, and must not
intersect the complement of the segment $S$ in the line through $a$ and $b$.\vadjust{\goodbreak}

%f6 ###
\begin{figure}

\includegraphics{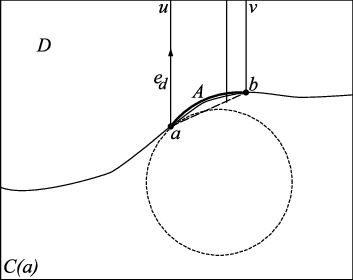}

\caption{Illustration of proof of Proposition
\protect\ref{proplipschitz}, part (\protect\hyperlink{itemlipschitzpart1-metric}{1}).
%KB new
The aspect ratio is not realistic---the height of $\cyl(a)$ must be at least
twice its horizontal diameter.}
\label{figshy2-2d-plane}
\end{figure}

The envelope of the boundaries of all such disks
%KB new
of radius $r\sin\alpha$ in $D\cap P$
%KB end new
is formed by the complement of the segment $S$ in the line through $a$
and $b$ together with
the minor arc $A$ of the circle of radius $r\sin\alpha$ running through
$a$ and $b$.
We can use $A$ to generate a short path between $a$ and $b$ in
$\overline{D}$ as follows.
If $A$ does not intersect one of the rays in \eqref{eqrays} then $A$
itself suffices; otherwise a still shorter path may be formed which
lies wholly in $\overline{D}$ by making a
short-cut using the relevant ray.
In any case a small perturbation of $A$ or the short-cut version, using
the vector $e_d$, will
provide a path in $D$ from $a$ to $b$ of length less than the length of
$A$ plus an arbitrarily small increment.
Calculation of the length of the minor arc $A$ now leads to the desired
bounds on $\intrinsicdist(a,b)$
as given in \eqref{eqlipschitzpart1-metric}.
\end{pf*}

\begin{pf*}{Proof of part (\protect\hyperlink{itemlipschitzpart2-differentiable}{2})}
Consider points $a,b$ and $c$ in $\ol D$, lying in this order
along an intrinsic geodesic
$\Gamma$
in $\ol D$.
% WSK 7 July 2011 we need $\Gamma$ to be minimal!
We will need the geodesic $\Gamma$ to be minimal. This is immediate in
case $\kappa=0$; in
the case $\kappa>0$ it follows if we require that the length of $\Gamma
$ is strictly less than $\pi/\sqrt{\kappa}$.
For some positive $t<\min\{\delta/(4\lambda),2r\sin\alpha\}$, suppose
that the intrinsic
distances between $a$ and $b$
and
between $b$ and $c$
are both equal to $t$.
Since $\Gamma$ is a minimal geodesic, the intrinsic distance
between $a$ and $c$ must be $2t$.
Let $\rho_1 = |a-b|$, $\rho_2 = |b-c|$ and $\rho_3
= |a-c|$ be the Euclidean distances between
these three
pairs of
points and
let $\pi- \theta$ be the
interior
angle at $b$ in the Euclidean triangle $abc$.
By the cosine formula,
\[
\cos\theta
= - \cos(\pi-\theta)
= - \frac{ \rho_1^2 + \rho_2^2 - \rho_3^2}{2 \rho_1 \rho_2}
= \frac{\rho_3^2 - \rho_1^2 - \rho_2^2}{2 \rho_1 \rho_2} .
\]
The upper bound on $t$ means
we can apply \eqref{eqlipschitzpart1-metric} to
the intrinsic and Euclidean distances between
$a$, $b$
and $c$.
Hence
$\rho_1 \leq t$, $\rho_2 \leq t$ and
\[
\rho_3 \geq 2 r \sin\alpha\sin\biggl( \frac{2 t}{2 r \sin\alpha}
\biggr)
\geq
2t \biggl( 1 - \frac{1}{6} \biggl( \frac{t}{r\sin\alpha}\biggr)^2 \biggr) ,
\]
where the last step uses $\sin\alpha\geq\alpha- \alpha^3/6$
if $\alpha\geq0$.
Together with the cosine formula, these
bounds for $\rho_1,\rho_2 $ and $\rho_3$
yield
\begin{eqnarray*}
\cos\theta &\geq&
\frac{(2t ( 1 - ({1}/{6}) ( {t}/{(r\sin\alpha)})^2 ))^2-\rho
_1^2-\rho_2^2}{2 \rho_1 \rho_2}
\\
&\geq&
\frac{(2t ( 1 - ({1}/{6} )( {t}/{(r\sin\alpha)})^2 ))^2-2t^2}{2 t^2}
=
2 \biggl( 1 - \frac{1}{6}
\biggl( \frac{t}{r\sin\alpha}\biggr)^2 \biggr)^2 -1 ,
\end{eqnarray*}
hence
\[
\cos\frac{\theta}{2} \geq 1 - \frac{1}{6}\biggl(\frac{t}{r\sin\alpha}\biggr)^2 .
\]
Considering
$t<\min\{\delta/(4\lambda),2r\sin\alpha\}$,
it follows by calculus that there exists a~$c(t)$ tending to zero with
$t$ such that
%
%e10 ###
\begin{equation}\label{eqangle-bound}
\theta \leq
\frac{4}{\sqrt{3}}
\frac{t}{2r\sin\alpha}\biggl(1+c(t)\frac{t}{2r\sin\alpha}\biggr) .
\end{equation}
%
% M1 - regarding the last paragraph: There are two places that you
%suggested insertions in response to questions of mine
% on Tuesday - the factor of 2 of d in one formula and the -1 term on
%the next line. I thought it better to let you
% insert these. Today, I questioned the exact form of (4) at the end of
%the paragraph. end M1
% WSK1/MB1: The inequality for theta here is correct though not as
%strong as it could be
% (could be linear + cubed correction in rho). This was checked with
%Mma.

Suppose now that the intrinsic geodesic $\Gamma$ has total length $K$.
For any positive integer $m>K/\min\{\delta/(4\lambda),2r\sin\alpha\}$,
let $a_0=a$, $a_1, \ldots, a_{m-1}$, $a_m=b$ be $m+1$ points
equally spaced along the geodesic, so that $\intrinsicdist
(a_{j-1},a_j)=t$ for $j=1, \ldots, m$.
Define $g_m\dvtx [0,K]\to\Reals^d$ to be the piecewise-linear curve interpolating
$g_m(jK/m)=a_j$
for $j=0, \ldots, m$.
By \eqref{eqangle-bound}, all the angles between successive
line-segments of the trajectory of $g_m$ are bounded above by
\[
\frac{4}{\sqrt{3}}
\frac{t}{2r\sin\alpha}\biggl(1+c(t)\frac{t}{2r\sin\alpha}\biggr) .
\]
Define the \textit{directional unit vector field} of the curve $g_m$ by
$\omega_m(s)=g_m^\prime(s)/\break |g_m^\prime(s)|$
for $s$ where $g_m(s)$ is linear, and extend to all $s$ using
left-limits for $s>0$ and the right-limit for $s=0$.
Then, by the triangle inequality,
\[
|\omega_m(s_2)-\omega_m(s_1)| \leq
\frac{4}{\sqrt{3}}
\frac{t}{2r\sin\alpha}\biggl(1+c(t)\frac{t}{2r\sin\alpha}\biggr)\biggl(\frac
{|s_2-s_1|}{t}+1\biggr) .
\]
From \eqref{eqlipschitzpart1-metric},
\[
\frac{2r\sin\alpha}{t} \sin\biggl(\frac{t}{2r\sin\alpha}\biggr) \leq |g_m^\prime
(s)| \leq 1 ;
\]
hence we obtain the inequality
%
%e12 ###
%e11 ###
\begin{eqnarray}\label{eqbdd-var-bound}
|g_m^\prime(s_2)-g_m^\prime(s_1)|
&\leq&|\omega_m(s_2)-\omega_m(s_1)|
+ |\omega_m(s_2)-g_m^\prime(s_2)|\nonumber\\
&&{}
+ |g_m^\prime(s_1)-\omega_m(s_1)|
\nonumber
\\[-8pt]
\\[-8pt]
\nonumber
&\leq&\frac{4}{\sqrt{3}}
\frac{t}{2r\sin\alpha}\biggl(1+c(t)\frac{t}{2r\sin\alpha}\biggr)\biggl(\frac{|s_2-s_1|}{t}+1\biggr)
\\
&&{}+
2\biggl(1- \frac{2r\sin\alpha}{t} \sin\biggl(\frac{t}{2r\sin\alpha}\biggr)\biggr)
,\nonumber
\end{eqnarray}
from which there follows a uniform bound on the absolute variation of
the $g_m^\prime$ functions.
Thus we can apply Helly's selection theorem to deduce that $g_m^\prime$
will converge along a subsequence, both pointwise and locally in $L^1$,
to a continuous limit $h$. It is immediate that $g_m$ converges
uniformly to $\Gamma$, and $\Gamma$ must be almost everywhere
differentiable with limit $h=\Gamma^\prime$. Moreover, from~\eqref
{eqbdd-var-bound} [and bearing in mind that $c(t)\to0$ with $t$],
we may deduce that the derivative $\Gamma^\prime$ is Lipschitz with constant
\[
\frac{4}{\sqrt{3}}\frac{1}{2r\sin\alpha} ,
\]
and indeed that $\Gamma$ is continuously differentiable.
\end{pf*}

\begin{pf*}{Proof of part (\protect\hyperlink{itemlipschitzpart3-chi}{3})}
As noted above, the
%KB $\CAT0$
$\CAT\kappa$
property of $\overline{D}$ implies that all geodesics
%KB in $\overline{D}$
between points $x$ and $y$ satisfying $0<\intrinsicdist(x,y) < \pi/\sqrt
{\kappa}$
are unique and minimal. Consider
%KB $(x,y),(x_n,y_n)\in(\overline{D}\times
%KB \overline{D})\setminus\{(x,x):x\in\overline{D}\}$
$(x,y),(x_n,y_n)\in\{(v,z) \in\overline{D}\times
\overline{D}\dvtx  0<\break \intrinsicdist(v,z) < \pi/\sqrt{\kappa} \}$
with $x_n\to x$ and
$y_n\to y$ in the Euclidean metric; taking subsequences we may suppose that
$\chi(x_n,y_n)$ converges to a limit. Part (\hyperlink{itemlipschitzpart2-differentiable}{2}) of the lemma
establishes the uniform
Lipschitz property of the direction fields of all minimal geodesics in
$\overline{D}$ so, by the Arzela--Ascoli theorem, we can find a subsequence
$(x_{n_k}, y_{n_k})$ such that the
geodesics from $x_{n_k}$ to $y_{n_k}$ must converge to a curve from
$x$ to $y$ whose direction field is the limit of the direction fields of
these minimal geodesics; hence its direction at $x$ must be
$\lim_k\chi(x_{n_k},y_{n_k})$. By minimality of the geodesics from
$x_{n_k}$ to $y_{n_k}$
and taking limits,
the length of the limiting curve can be no greater than that of the unique
minimal geodesic from~$x$ to~$y$; therefore the limiting curve must
also be
a minimal geodesic from~$x$ to~$y$. By the
%KB $\CAT0$
$\CAT\kappa$
property, the two
minimal geodesics from~$x$ to~$y$ must therefore be equal, and
therefore it
follows that $\lim_k\chi(x_{n_k},y_{n_k})=\chi(x,y)$. It
follows that any subsequence of $(x_n,y_n)\to(x,y)$ (convergence in Euclidean
metric) must possess a further subsequence for which
$\lim_k\chi(x_{n_k},y_{n_k})=\chi(x,y)$, and therefore
$\lim_n\chi(x_n,y_n)=\chi(x,y)$ must hold.
This establishes continuity of $\chi$ with reference to the Euclidean metric.
\end{pf*}

\begin{rem*}
Part (\hyperlink{itemlipschitzpart1-metric}{1}) of Proposition~\ref
{proplipschitz} may be used to show that $\chi(a,b)$ is H{\"o}lder($\frac{1}{2}$)
in its second argument $b$ when $a$ and $b$ are well-separated. We omit
this argument, as the result is not used in this paper.
\end{rem*}

\begin{rem*}\label{remtaylor}
Setting $\rho=|a-b|$ and $t=\intrinsicdist(a,b)$,
Inequality \eqref{eqlipschitzpart1-metric} can be rewritten as
\[
\frac{2r\sin\alpha}{d}\sin\biggl(\frac{t}{2r\sin\alpha}\biggr)
\leq
\frac{\rho}{t}
\leq
1 .
\]
The following is a trivial but useful consequence of the above estimates:
for some $c_1>0$, depending on $D$, and all $a,b\in\ol D$, with $|x-y|
\leq c_1$,
%
%e13 ###
\begin{equation}\label{m232}
|a-b| \leq\intrinsicdist{(a,b)} \leq2 |a-b|.\vadjust{\goodbreak}
\end{equation}
Moreover, since $\sin\phi\geq\phi-\frac{1}{6}\phi^3$ if $\phi\geq0$,
\[
1-\frac{1}{6}\frac{t^2}{4r^2\sin^2\alpha}
\leq
\frac{\rho}{t}
\leq
1 .
\]
The last inequality and \eqref{m232} imply that for some $c_2,c_3<\infty
$, depending on $\delta,r$ and~$\alpha$,
and for $\rho<\min\{\delta,2r\sin\alpha\}$,
%
%e14 ###
\begin{equation}\label{eqtaylor2}
1
\leq
\frac{t}{\rho}
\leq
1+c_2 t^2 \leq
1+c_3 \rho^2 .
\end{equation}
\end{rem*}

Proposition~\ref{proplipschitz} makes it possible to quantify the
extent to which short intrinsic geodesics may be approximated by
Euclidean segments.
\begin{cor}\label{corstiff}
% WSK 7 July 2011: have to keep propagating the substitution $\delta\to
Suppose the assumptions on Proposition~\ref{proplipschitz} hold, and
that~$\Gamma$ is a unit-speed intrinsic geodesic
with
intrinsic length $t<\min\{\delta/(4\lambda),\break 2r\sin\alpha\}$. Then
\[%\label{eqstiff}
|\Gamma(t)-\Gamma(0)-\Gamma^\prime(0) t|
\leq
\frac{4}{3}\frac{t^2}{2r\sin\alpha} .
\]
\end{cor}

\begin{pf}
Set $\rho=|\Gamma(t)-\Gamma(0)|$ equal to the Euclidean distance
between the two end-points of $\Gamma$;
then $\rho$ is bounded above by the intrinsic length $t$.
Let $\theta$ be the angle between $\Gamma^\prime(0)$ and $\Gamma
(t)-\Gamma(0)$.

Proposition~\ref{proplipschitz}(\hyperlink{itemlipschitzpart2-differentiable}{2}) tells us that $\Gamma^\prime$
is Lipschitz with constant $\frac{4}{\sqrt{3}}\frac{1}{2r\sin\alpha}$.
% Hence the angular distance
% between the unit vectors $\Gamma^\prime(0)$ and $\Gamma^\prime(s)$
% is bounded above by $\frac{4}{\sqrt{3}}\frac{s}{2r\sin\alpha}$
% when $s\geq0$.
% Using $\cos u\geq1-\frac{1}{2}u^2$, it follows that
% \begin{equation}\label{eqstiff-angle}
% \langle\Gamma^\prime(s), \Gamma^\prime(0) \rangle
% \geq
% 1 -\frac{8}{3} \frac{s^2}{4r^2\sin^2\alpha}
% \end{equation}
Hence
%
%e15 ###
\begin{equation}\label{eqstiff-angle}
\langle\Gamma^\prime(s), \Gamma^\prime(0) \rangle
=
1 - \frac{1}{2}|\Gamma^\prime(s)-\Gamma^\prime(0)|^2
\geq
1 -\frac{1}{2}\biggl(\frac{4}{\sqrt{3}}\frac{s}{2r\sin\alpha}\biggr)^2
\end{equation}
and this integrates to
\[
\langle\Gamma(t)-\Gamma(0), \Gamma^\prime(0) \rangle
\geq
\biggl(1 -\frac{8}{9} \frac{t^2}{4r^2\sin^2\alpha}\biggr) t .
\]
Consequently
\begin{eqnarray*}
|\Gamma(t)-\Gamma(0)-\Gamma^\prime(0) t|^2
&=&
|\Gamma(t)-\Gamma(0)|^2 + |\Gamma^\prime(0) t|^2
-2\langle\Gamma(t)-\Gamma(0),\Gamma^\prime(0) t\rangle
\\
&\leq&
\rho^2 + t^2 - 2 \biggl(1 -\frac{8}{9} \frac{t^2}{4r^2\sin^2\alpha}\biggr) t^2
\leq
\frac{16}{9} \frac{t^4}{4r^2\sin^2\alpha} .
\end{eqnarray*}
The result follows by taking square roots.
\end{pf}

At this point we revert to considering $\CAT0$ spaces only, since
generalization of the following proofs to the $\CAT\kappa$ case would
extend the exposition.
We recall Gauss' lemma from Riemannian geometry, that the exponential
map is a radial isometry.
\citeauthor{CheegerEbin-1975} [(\citeyear{CheegerEbin-1975}), Chapter 1, Section 2] observe that,
for smooth Riemannian manifolds, it is equivalent to the assertion that
the Riemannian distance
$\intrinsicdist(x,y)$ is continuously differentiable
in $x$ when $x\neq y$ and $y$ does not lie in the cut-locus of $x$,
with the gradient being given by the tangent of the geodesic running
from $y$ to $x$.
Proposition~\ref{proplipschitz} and Corollary~\ref{corstiff} can be
used to prove the following Gauss lemma for $\CAT0$
domains with sufficient boundary regularity. Here, $\grad
_x\intrinsicdist(x,y)$ refers to the Euclidean gradient with
respect to~$x$, with $\grad\intrinsicdist(x,y)$ being the gradient with
respect to both variables.

%f7 ###
\begin{figure}

\includegraphics{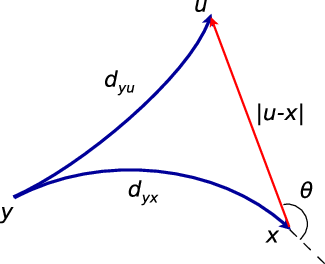}

\caption{Illustration of the configuration of the triangle referred to
in the statement of Proposition~\protect\ref{propgauss}.
The sides running from $y$ to $u$ and from $y$ to $x$ (and of intrinsic
lengths $d_{yu}$ and $d_{yx}$, resp.)
are intrinsic geodesics. The side running from $x$ to $u$
is a Euclidean segment.}
\label{figshy2-gauss1}
\end{figure}

Note also that a consequence of Proposition~\ref{proplipschitz} is that
intrinsic geodesics have continuously varying directions, and therefore
that it makes sense to speak of the angle between a geodesic and a
Euclidean segment.

\begin{prop}\label{propgauss}
Suppose
% (as in Proposition~\ref{proplipschitz})
that $D$ is
a $\CAT0$
domain,
bounded in the Euclidean metric,
satisfying a uniform exterior sphere
condition based on radius $r>0$, and a uniform interior cone
condition based on radius $\delta>0$ and angle
$\alpha\in(0,\pi/2]$.
For every $c_1 >0$, there exist $c_2,c_3<\infty$ such that, if $x,y \in
\ol D$ with
$\intrinsicdist(x,y)\geq c_1 $ and $\sqrt{|u-x|} \lor|u-x| \leq c_2$, then
%
%e16 ###
\begin{equation}\label{eqgauss2}
\bigl|\intrinsicdist(u,y)-\bigl(\intrinsicdist(x,y)+|u-x|\cos\theta\bigr)\bigr|
\leq c_3 |u-x|^{3/2},
\end{equation}
where $\theta$ is the angle between the geodesic from $y$ to $x$
and the Euclidean segment from $x$ to $u$ that is exterior to the direction
from $y$ to $x$ (see Figure~\ref{figshy2-gauss1}).
Consequently, if $x$, $y\in\ol D$ with $x\neq y$, then
%
%e17 ###
\begin{equation}\label{eqgauss1}
\grad_x\intrinsicdist(x,y) = -\chi(x,y).
\end{equation}
Moreover,
%
%e18 ###
\begin{equation}\label{eqgauss3}
\grad\intrinsicdist(x,y) = (-\chi(x,y), -\chi(y,x)).
\end{equation}
\end{prop}
%
%WSK3: start

Note that \citet{Bieske-2010} establishes a similar result for
Carnot--Carath\'eo\-dory spaces. In both cases, the relevant distance function
satisfies an eikonal equation.
%WSK3: end

\begin{pf*}{Proof of Proposition \protect\ref{propgauss}}
In order to demonstrate \eqref{eqgauss2}, we establish upper and lower
bounds on the difference
%
%e19 ###
\begin{equation}\label{eqdistintr}
\intrinsicdist(y,u)-\intrinsicdist(y,x)
\end{equation}
when $u$ is close to $x$. We abbreviate, setting $d_{yx}=\intrinsicdist
(y,x)$, {etc}.

Let $\theta'$ be the exterior angle between the geodesic from $x$ to $y$
and the geodesic from $x$ to $u$.
By the $\CAT0$ property, the Euclidean triangle with the same side
lengths as a triangle
in the intrinsic metric has larger interior angles and therefore
smaller exterior angles. [The elementary argument for this is given in
\citet{BridsonHaefliger-1999}, Chapter II.1, Proposition~1.7(4).]
Thus if $\theta''$ is the exterior
angle of the comparison triangle for $x,y$ and $u$
corresponding to the exterior angle $\theta'$, then $\theta''\leq\theta
'$, and so
\begin{eqnarray*}
d_{yu} &=& \sqrt{d_{yx}^2+d_{xu}^2+2d_{yx}d_{xu}\cos\theta''}
\geq
\sqrt{d_{yx}^2+d_{xu}^2+2d_{yx}d_{xu}\cos\theta'}\\
&\geq&
\sqrt{d_{yx}^2+d_{xu}^2\cos^2\theta'
+2d_{yx}d_{xu}\cos\theta'}
= d_{yx}+d_{xu}\cos\theta' \\
& \geq& d_{yx}+|u-x|\cos\theta' .
\end{eqnarray*}
Corollary~\ref{corstiff}
implies that
$|\theta- \theta'| \leq c_4 d_{xu}$ for small $d_{xu}$. Hence,
$|\theta- \theta'| \leq c_5 |u-x|$ and
$|\cos\theta- \cos\theta'| \leq c_5 |u-x|$.
We obtain for $|u-x|\le c_2$, for some $c_2>0$,
%
%e20 ###
\begin{eqnarray} \label{m234}
\qquad d_{yu} &\geq &d_{yx}+|u-x|\cos\theta'
\geq d_{yx}+|u-x|\cos\theta
-|u-x||\cos\theta- \cos\theta'|
\nonumber
\\[-8pt]
\\[-8pt]
\nonumber
&\geq& d_{yx}+|u-x|\cos\theta
-c_5|u-x|^2 .
\end{eqnarray}
This provides a lower bound on (\ref{eqdistintr}) and a bound for one
direction of \eqref{eqgauss2}.

We now establish an upper bound on (\ref{eqdistintr}). Fix a point $w$
on the intrinsic geodesic from $y$ to $x$.
Then $d_{yx}=d_{yw}+d_{wx}$ and $d_{yu}\leq d_{yw}+d_{wu}$.
We shall require $w$ to be close to $x$, but not as close as $u$, with
$|w-x|=\sqrt{|u-x|}$ being assumed.

Because $w$ is close to $x$ and thus also close to $u$, we may replace
the intrinsic geodesics from $w$ to $u$ and
from $w$ to $x$ by Euclidean segments, without greatly altering lengths
and segments. Let $\theta^*$ be the exterior
angle at $x$ for the Euclidean triangle with sides $d_{x,w}, d_{x,u}$
and $d_{u,w}$.
From \eqref{eqtaylor2},
%
%e22 ###
%e21 ###
\begin{eqnarray}
|u-w| &\leq& d_{wu} = \frac{d_{wu}}{|u-w|}|u-w| \leq
(1+ c_6 |u-w|^2 ) |u-w|,\\
|x-w| &\leq& d_{wx} = \frac{d_{wx}}{|x-w|}|x-w| \leq
(1+ c_6 |x-w|^2 ) |x-w|,
\end{eqnarray}
when $|u-w|$, $|x-w| < 2 r \sin\alpha$.\vadjust{\goodbreak}

As before, by Corollary~\ref{corstiff},
$|\theta- \theta^*| \leq c_7 d_{xw}$ for small $d_{xw}$. Hence,
$|\theta- \theta^*| \leq c_8 |w-x|$ and
%
%e23 ###
\begin{equation}\label{m233}
|\cos\theta- \cos\theta^*| \leq c_8 |w-x|.
\end{equation}
These computations allow use to establish an upper bound for
$\intrinsicdist(y,u)-\intrinsicdist(y,x)$.
First note that
\begin{eqnarray*}
d_{yu} &\leq& d_{yw}+d_{wu} \leq d_{yw} + (1+ c_6 |u-w|^2) |u-w|\\
&\leq&
d_{yw} + |u-w| + c_6 |u-w|^3\\
&\leq&
d_{yw} + |u-w| + c_6 (|u-x|+|w-x|)^3.
\end{eqnarray*}
Now apply the cosine formula to control $|u-w|$, using \eqref{m233}:
\begin{eqnarray*}
|u-w| &= &\sqrt{|w-x|^2+|u-x|^2+2|w-x||u-x|\cos\theta^*}\\
&=&
\bigl((|w-x|+|u-x|\cos\theta)^2+|u-x|^2\sin^2\theta\\
&&\hspace*{40pt}{}+2|w-x||u-x|(\cos
\theta^*-\cos\theta)\bigr)^{1/2}\\
&\leq&
|w-x|+|u-x|\cos\theta\\
&&{}+\frac{1}{2}\frac{|u-x|^2\sin^2\theta}{|w-x|+|u-x|\cos\theta}
+\frac{|w-x||u-x|(\cos\theta^*-\cos\theta)}{|w-x|+|u-x|\cos\theta}\\
&\leq&
|w-x|+|u-x|\cos\theta
+\frac{1}{2}\frac{|u-x|^2\sin^2\theta}{|w-x|-|u-x|}
\\
&&{}+\frac{|w-x||u-x|}{|w-x|-|u-x|}
c_8 |w-x|.
\end{eqnarray*}
If we take
$|w-x|=\sqrt{|u-x|}$, with $|u-x|< c_9$ for a suitably small $c_9>0$, then
\[
|u-w| \leq |w-x| + |u-x|\cos\theta+ c_{10} |u-x|^{3/2}.
\]
Combining these bounds implies
\begin{eqnarray*}
d_{yu} &\leq&
d_{yw} + |u-w| + c_6 (|u-x|+|w-x|)^3 \\
&\leq&
d_{yw} + |w-x| + |u-x|\cos\theta+ c_{10} |u-x|^{3/2} + c_6
(|u-x|+|w-x|)^3 \\
&\leq&
d_{yw} + d_{w,x} + |u-x|\cos\theta+ c_{11} |u-x|^{3/2}\\
&=&
d_{yx} + |u-x|\cos\theta+ c_{11} |u-x|^{3/2},
\end{eqnarray*}
which provides an upper bound on (\ref{eqdistintr}).
It follows from the above inequality and~\eqref{m234} that
\[
\bigl|d_{yu}-(d_{yx} + |u-x|\cos\theta)\bigr| \leq c_3 |u-x|^{3/2} ,
\]
which yields the bound in \eqref{eqgauss2}. The formula in \eqref
{eqgauss1}, for the gradient of the
intrinsic distance $\intrinsicdist(x,y)$ with respect to $x$, follows
immediately.

We still need to demonstrate the formula in (\ref{eqgauss3}).
Let $c_1, c_2$ and $c_3$ be as in the statement of \eqref{eqgauss2}.
Fix $x,y \in\ol D$ and suppose that $\intrinsicdist{(x,y)} \geq2
c_1$. Suppose that $u,v \in\ol D$,
$\sqrt{|u-x|} \lor|u-x| \leq c_2 \land c_1/4$ and
$\sqrt{|v-y|} \lor|v-y| \leq c_2\land c_1/4$.
Let $\theta_x$ be the exterior angle between the geodesic from $x$ to $y$
and the Euclidean segment from $x$ to $u$.
Similarly, let $\theta_y$ be the exterior angle between the geodesic
from $y$ to $x$
and the Euclidean segment from $y$ to $v$.
Also, let $\theta'$ be the exterior angle between the geodesic from $y$
to $u$
and the Euclidean segment from $y$ to $v$.
Then by the above reasoning
%
%e24 ###
\begin{equation}\label{m236}
\bigl|\intrinsicdist(u,y)-\bigl(\intrinsicdist(x,y)+|u-x|\cos\theta_x\bigr)\bigr|
\leq c_3 |u-x|^{3/2}
\end{equation}
and
%
%e25 ###
\begin{equation}\label{m235}
\bigl|\intrinsicdist(v,u)-\bigl(\intrinsicdist(y,u)+|v-y|\cos\theta'\bigr)\bigr|
\leq c_3 |v-y|^{3/2}.
\end{equation}

Recall that the Euclidean triangle with the same side lengths as a
triangle in intrinsic metric
has larger interior angles. Using a triangle inequality for angles,
$|\theta_y - \theta'| $ is less than the angle at the
vertex corresponding to $y$ in the Euclidean triangle with sides
$d_{xy}, d_{xu}$ and $d_{yu}$. It follows
that $|\theta_y - \theta'| \leq c_{12} d_{xu}/d_{xy} \leq c_{13} |u-x|$
and therefore
$|\cos\theta_y - \cos\theta'| \leq c_{13} |u-x|$.
This and \eqref{m235} yield
%
%e26 ###
\begin{eqnarray}\label{m237}
&&\bigl|\intrinsicdist(v,u)-\bigl(\intrinsicdist(y,u)+|v-y|\cos\theta_y\bigr)\bigr|
\nonumber
\\[-8pt]
\\[-8pt]
\nonumber
&&\qquad\leq c_3 |v-y|^{3/2} + c_{13}|u-x||v-y|.
\end{eqnarray}
The triangle inequality applied to the left-hand sides of \eqref{m236}
and \eqref{m237} implies that
\begin{eqnarray*}%\label{m237}
&& \bigl|\intrinsicdist(v,u)-\bigl(\intrinsicdist(x,y)+|v-y|\cos\theta_y +
|u-x|\cos\theta_x\bigr)\bigr| \\
&&\quad \leq c_3 |u-x|^{3/2} + c_3 |v-y|^{3/2} + c_{13}|u-x||v-y|.
\end{eqnarray*}
Consequently, $\grad\intrinsicdist(x,y)$ exists when $\intrinsicdist
(x,y)$ is viewed
as a function of $(x,y)\in(\overline{D}\times\overline{D})\setminus\{
(u,u)\dvtx u\in\overline{D}\}$ and is given by
\[
\grad\intrinsicdist(x,y) = (-\chi(x,y), -\chi(y,x)) .
\]
\upqed\end{pf*}

In Proposition~\ref{propwell-posed}, we consider solutions of the
differential equation $\d x=\chi(x,y)\,\d t$
for pursuit and evasion.
Proposition~\ref{proplipschitz} established partial regularity
for $\chi(x,y)$, which does not automatically
guarantee well-posedness of solutions (as defined below).
However, the $\CAT0$ property, together with boundary regularity,
will imply well-posedness, even for some discontinuous driving paths $y$.

Suppose that $y(t)$, $t\in[0,T_1]$, is cadlag, of bounded variation on
finite intervals, and takes values in $\ol D$.
We will say that $x(t)$, $t\in[0, T_1]$, is a weak solution to
$ \d x=\chi(x,y)\,\d t$ if $x(t) = x(0) + \int_0^t
\chi(x(s),y(s))\,\d s$ for all $t\in[0,T_1]$.

\begin{prop}\label{propwell-posed}
Let $D$ be a $\CAT0$ domain satisfying uniform exterior sphere and
interior cone conditions.
For distinct $x,y\in\ol D$, let $\chi(x,y)$ be the unit tangent vector
at $x$ of the geodesic
from $x$ to $y$. We consider the differential equation
%
%e27 ###
\begin{equation}\label{m2315}
\d x = \chi(x,y)\,\d t
\end{equation}
defined in the weak sense for absolutely continuous paths $\{x(t)\dvtx t \geq
0\}$ in $\ol D$,
driven by paths $\{y(t)\dvtx t \geq0\}$,
up until the first time that $x$ and $y$ are equal.
The problem is well-posed, in the sense that solutions $x$ exist,
are uniquely determined by initial values $x(0)$,
and depend continuously on the initial value $x(0)$ and the
driving process $\{y(t)\dvtx t \geq0\}$ (using the uniform distance metric
in both cases).
\end{prop}

\begin{pf}
The argument is based on the simpler case when the path $y$ is constant
in time,
which we for the moment assume. In this case, existence
follows directly from the existence of intrinsic geodesics in $\CAT0$ domains.
To show uniqueness, note that, for two solutions $x$ and $\widetilde
{x}$ of (\ref{m2315}),
since $x$
and $\widetilde{x}$ are absolutely continuous and satisfy the
differential equation weakly,
for almost all $s$, the time-derivatives of $x(s)$ and $\widetilde
{x}(s)$ must exist
and be given by $\chi(x(s),y)$ and $\chi(\widetilde{x}(s),y)$.
Exploiting the differentiability of the intrinsic distance given by
Proposition~\ref{propgauss},
for $x(s)\neq y$ and $\widetilde{x}(s)\neq y$, one has
%
%e28 ###
\begin{equation}\label{m2311}
\qquad\biggl[\frac{\d}{\d t}\intrinsicdist\bigl(x(s+t),\widetilde{x}(s+t)\bigr)\biggr]_{t=0} =
\biggl[\frac{\d}{\d t}\intrinsicdist\bigl(\Gamma^{(s)}(t),\widetilde{\Gamma
}^{(s)}(t)\bigr)\biggr]_{t=0} ,
\end{equation}
where $\Gamma^{(s)}$, $\widetilde{\Gamma}^{(s)}$ are unit-speed
geodesics running
from $x(s)$, $\widetilde{x}(s)$ to $y$.
We will show that
%
%e29 ###
\begin{equation}\label{m2310}
\biggl[\frac{\d}{\d t}\intrinsicdist\bigl(\Gamma^{(s)}(t),\widetilde{\Gamma
}^{(s)}(t)\bigr)\biggr]_{t=0}
\leq0 .
\end{equation}
Consider a Euclidean triangle $abc$ with side lengths satisfying
$|ab| = \break \intrinsicdist({\Gamma}^{(s)}(0), y)$,
$|cb|=\intrinsicdist(\widetilde{\Gamma}^{(s)}(0), y)$ and
$|bc|=\intrinsicdist({\Gamma}^{(s)}(0), \widetilde{\Gamma}^{(s)}(0))$.
Let $z(t)\in ab$ be a point such that $|z(t) - a| = t$, and let
$\wt z(t)\in bc$ be a point such that $|\wt z(t) - c| = t$.
Then Definition~\ref{m221}
implies that
\[
\intrinsicdist\bigl({\Gamma}^{(s)}(t), \widetilde{\Gamma}^{(s)}(t)\bigr)
\leq |z(t) - \wt z(t)|
\leq |z(0) - \wt z(0)|
= \intrinsicdist\bigl({\Gamma}^{(s)}(0), \widetilde{\Gamma}^{(s)}(0)\bigr).
\]
This implies \eqref{m2310}. It follows that the derivative on the left-hand side of \eqref{m2311} is nonpositive;
therefore $x=\widetilde{x}$ if $x(0)=\widetilde{x}(0)$, and so
uniqueness holds.

By considering the behaviour over disjoint time intervals, existence
and uniqueness follow for the case when $y$ is piecewise-constant,
in which case the solution curve $x$ is piecewise-geodesic.

We will establish continuous dependence on the initial position $x(0)$
and the driving process $y$,
when $y$ is piecewise constant.
Suppose that\vadjust{\goodbreak} $y$, $\widetilde{y}$ are
two piecewise-constant paths in $\ol D$, and $x$, $\widetilde{x}$ solve
\[
\d x = \chi(x,y)\,\d t ,\qquad
\d\widetilde{x} = \chi(\widetilde{x},\widetilde{y})\,\d t
\]
for prescribed initial positions $x(0) \neq y(0)$ and $\widetilde{x}(0)
\neq\widetilde{y}(0)$.
The solutions $x$, $\widetilde{x}$ satisfy the differential equations
weakly, and therefore,
for almost all $s$, the time-derivatives of $x(s)$ and $\widetilde
{x}(s)$ must exist
and are given by $\chi(x(s),y(s))$ and $\chi(\widetilde{x}(s),\widetilde
{y}(s))$.
% =================
Arguing as before,
for $x(s)\neq y(s)$ and $\widetilde{x}(s)\neq\widetilde{y}(s)$, we may
construct a $\CAT0$ comparison
for the two triangles defined by (a) vertices $x(s)$, $\widetilde
{x}(s)$, $y(s)$ and (b)
vertices $\widetilde{x}(s)$, $\widetilde{y}(s)$, $y(s)$ (see Figure \ref
{figshy2-wp1}).
Using this comparison,
and continuing until either $x(t) = y(t)$ or $\widetilde{x}(t) =
\widetilde{y}(t)$,
the function
$\intrinsicdist(x,\widetilde{x})$ is dominated by its Euclidean
counterpart for a two-dimensional
quadrilateral which is based on a pair of opposing sides of lengths
$\intrinsicdist(x(s),\widetilde{x}(s))$
and $\intrinsicdist(y(s),\widetilde{y}(s))$.
% WSK change starts here.

%f8 ###
\begin{figure}

\includegraphics{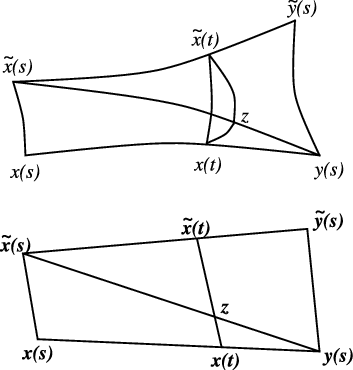}

\caption{Illustration of the $\CAT0$ comparison argument applied to the
triangles defined by vertices \textup{(a)} $x(s)$, $\widetilde{x}(s)$, $y(s)$
and \textup{(b)} $\widetilde{x}(s)$, $\widetilde{y}(s)$, $y(s)$.
The corresponding Euclidean triangles have vertices marked with
boldface symbols.}
\label{figshy2-wp1}\vspace*{-3pt}
\end{figure}

In detail, and using boldface symbols to indicate corresponding
Euclidean comparison points, we may argue as follows (see Figure \ref
{figshy2-wp1}).
Because side-lengths of comparison triangles agree,
\begin{eqnarray*}
\intrinsicdist(x(s), y(s)) &=& |\mathbf{x}(\mathbf{s})- \mathbf{y}(\mathbf{s})| ,\\
\intrinsicdist(\wt x(s), y(s)) &=& |\wt{\mathbf{x}}(\mathbf{s})- \mathbf{y}(\mathbf{s})| ,\\
\intrinsicdist(x(s), \wt x(s)) &=& |\mathbf{x}(\mathbf{s})- \wt{\mathbf{x}}(\mathbf{s})| ,\\
\intrinsicdist(\wt x(s), \wt y(s)) &=& |\wt{\mathbf{x}}(\mathbf{s})- \wt{\mathbf{y}}(\mathbf{s})|
,\\
\intrinsicdist(y(s),\wt y(s)) &= &|\mathbf{y}(\mathbf{s})- \wt{\mathbf{y}}(\mathbf{s})| .
\end{eqnarray*}
Locating $x(t)$ according to distance from $x(s)$ along the intrinsic
geodesic from $x(s)$ to $y(s)$,
and $\wt x(t)$ according to distance from $\wt x(s)$ along the\vadjust{\goodbreak}
intrinsic geodesic from $\wt x(s)$ to $\wt y(s)$
(and locating comparison Euclidean points in the corresponding way), we find
\begin{eqnarray*}
\intrinsicdist(x(s), x(t)) &=& |\mathbf{x}(\mathbf{s})- \mathbf{x}(\mathbf{t})| ,\\[-2pt]
\intrinsicdist(\wt x(s), \wt x(t)) &= &|\wt{\mathbf{x}}(\mathbf{s})- \wt{\mathbf{x}}(\mathbf{t})| .
\end{eqnarray*}
Now locate the Euclidean point $\mathbf{z}$ at the intersection of the
Euclidean line segments $\ol{\wt{\mathbf{x}}(\mathbf{s}), \mathbf{y}(\mathbf{s})}$ and $\ol{\mathbf{x}(\mathbf{t}),\wt{\mathbf{x}}(\mathbf{t})}$,
and locate $z$ on the intrinsic geodesic from $\wt x(s)$ to
$y(s)$ so that
\[
\intrinsicdist(\wt x(s), z) = |\wt{\mathbf{x}}(\mathbf{s})- \mathbf{z}| .
\]
Using comparison arguments and the nature of the Euclidean
parallelogram $\mathbf{x}(\mathbf{s})\wt{\mathbf{x}}(\mathbf{s})\wt{\mathbf{y}}(\mathbf{s}) \mathbf{y}(\mathbf{s})$, we then see that
\begin{eqnarray*}
\intrinsicdist(x(t), \wt x(t))
&\leq&\intrinsicdist(x(t), z) + \intrinsicdist(z, \wt x(t))\\[-2pt]
&\leq&|\mathbf{x}(\mathbf{t}) - \mathbf{z}| + |\mathbf{z}- \wt{\mathbf{x}}(\mathbf{t})|
= |\mathbf{x}(\mathbf{t}) - \wt{\mathbf{x}}(\mathbf{t})|
\\[-2pt]
&\leq&\max\{|\mathbf{x}(\mathbf{s}) - \wt{\mathbf{x}}(\mathbf{s})|,
|\mathbf{y}(\mathbf{s}) - \wt{\mathbf{y}}(\mathbf{s})|\}\\[-2pt]
&=& \max\{\intrinsicdist(x(s),\widetilde{x}(s)), \intrinsicdist
(y(s),\widetilde{y}(s))\} .
\end{eqnarray*}
This comparison can also be justified by use of Reshetnyak
majorization, however we have chosen to present an explicit elementary proof.
% WSK proposed change ends.

We now consider general $y$ in \eqref{m2315}.
There exist piecewise-constant functions~$y_n$ converging to $y$
uniformly on compact intervals;
let $x_n$ be the corresponding solutions to \eqref{m2315}, with $x_n(0)
= x(0)\in\ol D$.
If $|y_n(t) - y_m(t)| \leq c_1$ for $t\in[0, T]$, then $|x_n(t) -
x_m(t)| \leq c_1$
by the argument given above. Since the sequence $y_n$ is Cauchy in the
uniform norm
on $[0,T]$, so is the sequence $x_n$, which therefore converges to a
function $x$.

Recall that we are assuming $x(0) \ne y(0)$.
Choose fixed $\eps_1,\eps_2>0$ and let $T=\inf\{t>0\dvtx  |x(t-) -y(t-)|
\leq2\eps_1\}$.
By part (3) of Proposition~\ref{proplipschitz}, there exists \mbox{$\delta
_1>0$} such that,
if $|u_1 - u_2| \geq\eps_1$, $|v_1 - v_2| \geq\eps_1$, $|u_1 - v_1|
\leq\delta_1$ and
$|u_2 - v_2| \leq\delta_1$, then $|\chi(u_1,u_2) - \chi(v_1,v_2)| \leq
\eps_2/T$.
Suppose that $n$ is large enough so that $|y_n(t) - y(t)| \leq\delta
_1\land\eps_1$ for $t\in[0,T]$. Then
$|x_n(t) - x(t)| \leq\delta_1\land\eps_1$ for $t\in[0,T]$ and
$|\chi(x(t),y(t)) - \chi(x_n(t),y_n(t))| \leq\eps_2/T$.
We obtain, for $t\leq T$,
\begin{eqnarray*}
&&\biggl| x(t) -x(0) - \int_0^t \chi(x(s), y(s)) \,\d s\biggr|\\[-2pt]
&&\qquad \leq
|x(t) - x_n(t)| +
\biggl| x_n(t) - x(0) - \int_0^t \chi(x_n(s), y_n(s)) \,\d s \biggr|\\[-2pt]
&&\qquad\quad{} +
\biggl| \int_0^t \chi(x_n(s), y_n(s)) \,\d s - \int_0^t \chi(x(s), y(s)) \,\d s
\biggr|\\[-2pt]
&&\qquad \leq \delta_1 + 0 + \int_0^t |\chi(x_n(s), y_n(s)) - \chi(x(s),
y(s))| \,\d s \\[-2pt]
&&\qquad \leq \delta_1 + \eps_2.
\end{eqnarray*}
Since $\eps_1,\eps_2$ and $\delta_1$ can be chosen arbitrarily small,
we see that $x(t) = x(0) +\int_0^t \chi(x(s), y(s)) \,\d s$
for all $t < \inf\{t>0\dvtx  |x(t-) -y(t-)| =0\}$. Hence, $x$ is a solution
to \eqref{m2315}.

Uniqueness of the solution of \eqref{m2315}, for given $x(0)$ and
general $y$, follows by reasoning as in \eqref{m2311} and \eqref{m2310}.
The continuous dependence of solutions on $x(0)$ and $y$, for general
$y$, follows from the above estimates by approximating
$y$ by piecewise constant driving processes.
\end{pf}

% \section[CAT(0) and pursuit-evasion]{$\CAT0$ and pursuit-evasion}
%s3 ###
\section{\texorpdfstring{$\CAT0$ and pursuit--evasion.}{$\CAT0$ and pursuit--evasion}}\label{seccat0-and-pursuit-evasion}
We consider the Lion and Man problem in a bounded $\CAT0$ domain $D$
satisfying the uniform exterior sphere and
interior cone conditions.
\citet{AlexanderBishopGhrist-2006} showed that $\varepsilon/2$-capture,
for given $\varepsilon> 0$, must occur
for the discrete-time variant of this problem.
As we will see in Section~\ref{seccoupling-to-pursuit}, the Lion and
Man trajectories
$x$ and $y$
will be weak limits of couplings of reflected Brownian motions, with
drift and small noise, that arise from our
capture problem.

We therefore modify the \citet{AlexanderBishopGhrist-2006} argument
to apply to continuous time;
the modified argument also supplies an explicit upper bound on the
capture time.
We will only need to consider trajectories $x$ and $y$ that are
Lipschitz with constant $1$. Note that
Lipschitz trajectories are absolutely continuous, so that the
directions $\d x/\d t$ and $\d y/\d t$ are
defined for almost all times~$t$.

One can express the trajectories of Lion $x$ and Man $y$ as functions
of time
$t$ in the following differential form:
%
%e30 ###
\begin{eqnarray}\label{eqdeterministic-skorokhod}
\d x &=& \chi(x,y)\,\d t - \nu_x \,\d L^x ,
\nonumber
\\[-8pt]
\\[-8pt]
\nonumber
\d y &=& H \,\d t - \nu_y \,\d L^y .
\end{eqnarray}
Here, $H$ is assumed to be a pre-assigned, time-varying unit length
vector generating the
motion of the Man, $\chi(x,y)$ generates the motion of the Lion and is
defined as in
Proposition~\ref{proplipschitz}, for $x\neq y$, as the unit tangent at
$x$ for the corresponding
intrinsic geodesic, while
$\nu_x \in\NN_{x,r}$ and $\nu_y \in\NN_{y,r}$ (for $r>0$ satisfying
the exterior sphere condition of
$D$ as given in Definition~\ref{defexterior-sphere-condition})
determine the reflection off of the boundary $\partial D$.
The vector $H$ is assumed to be measurable
in $t$; on account of Proposition~\ref{proplipschitz}, $\chi$ is
continuous on $x \neq y$.
The terms $\nu_x \,\d L^x$,
respectively, $\nu_y \,\d L^y$, are differentials
arising from Skorokhod transformations and are differentials of
functions of bounded variation that increase only when $x$,
respectively, $y$, belong to $\partial D$, and
are then directed along an outward-pointing unit normal so as
to cancel exactly with the outward-pointing component of the
drifts $\chi\,\d t$, respectively, $H\,\d t$.

We note that Skorokhod
transformations are uniquely defined for a domain satisfying
uniform exterior sphere and interior cone conditions
[\citet{Saisho-1987}] [also compare earlier results of
\citet{LionsSznitman-1984}], and they then depend
continuously on the driving processes (using the uniform path metric).
In fact, by the definition of
$\chi$, $\nu_x \,\d L^x$ vanishes identically, while
$\nu_y \,\d L^y$ vanishes identically if $\langle H,
\nu\rangle<0$ whenever $y\in\partial D$.
In particular, Proposition~\ref{propwell-posed} applies and guarantees
the existence
of $x$ and its approximation by piecewise-geodesic paths for
$y$ determined by
%$x$.
$H$. [We include
both the Skorokhod transformation differentials in
\eqref{eqdeterministic-skorokhod} as they will both appear in
the stochastic version
in Section~\ref{seccoupling-to-pursuit}.]

We base our argument on \citeauthor{AlexanderBishopGhrist-2006} [(\citeyear{AlexanderBishopGhrist-2006}), Theorem~12].
The proof analyzes the greedy pursuit strategy
arising from the definition of the vector field $\chi$, with
the Lion always directing its motion along the intrinsic
geodesic from Lion to Man.
The $\CAT0$ property forces the
distance between Lion and Man to be
nonincreasing, and the Man must run directly away from the Lion
in order to prolong successful evasion. Since the domain is bounded,
this will, however, not be
achievable indefinitely.

In order to demonstrate the main result in this section, Proposition
\ref{propgreedy}, we will employ
the following lemma.
\begin{lem}\label{lemangle-pursuit}
Under the greedy pursuit strategy described above, in a\break $\CAT0$ domain
satisfying
uniform exterior sphere and interior cone conditions,
and at a time $t$ at which Lion and Man locations $x(t)$ and $y(t)$ are
differentiable in $t$,
\[
\frac{\d}{\d t}\intrinsicdist(x(t),y(t)) =
-\bigl(1-|y^\prime(t)|\cos\alpha(t)\bigr),
\]
where $\alpha(t)$ is the angle between the Man's velocity $y^\prime(t)$ and
the geodesic running from Lion to Man.
\end{lem}

\begin{pf}
This follows immediately from the generalization of Gauss' lemma to
such domains,
as was established in Proposition~\ref{propgauss}.
\end{pf}

%KB new
Alternatively, Lemma~\ref{lemangle-pursuit} follows directly from the
first variation formula in
$\CAT0$ spaces [\citet{BridsonHaefliger-1999}, page 185, \citet
{BuragoBuragoIvanov-2001}, Exercise~4.5.10].
%KB end new

\begin{prop}\label{propgreedy}
Suppose that $D$ is a bounded $\operatorname{CAT}(0)$ domain that satisfies a uniform
exterior sphere condition based on a radius $r>0$ and a uniform
interior cone condition based on a radius
$\delta> 0$ and angle $\alpha\in(0,\pi/2]$.
Under the greedy pursuit strategy described above,
there is a positive constant $t_c$ depending only on the
diameter of $D$ and $\varepsilon> 0$ [and not on $H$ in (\ref
{eqdeterministic-skorokhod})] such that the Lion will
come within distance $\varepsilon/2$ of the Man before time $t_c$,
regardless of their starting positions within $D$.
\end{prop}

\begin{rem*}\label{remepsilon2}
We use $\varepsilon/2$ here rather than $\varepsilon$, since a further
distance $\varepsilon/2$ will be required by
the stochastic part of the argument.
\end{rem*}

\begin{pf*}{Proof of Proposition \protect\ref{propgreedy}}
This proof follows \citet{AlexanderBishopGhrist-2006}, but is modified
(a) to account for the continuous time context and
(b) because we need to derive a specific upper bound $t_c$ on the time
of $\varepsilon/2$-capture.
Below, we abbreviate by setting $\mathcal{L}(t)=\intrinsicdist(x(t), y(t))$.

Let $\alpha$ be the angle defined in Lemma~\ref{lemangle-pursuit}.
Note that this is defined for almost all times $t$, since the paths
$x(t)$, $y(t)$ are Lipschitz and
are therefore differentiable for almost all $t$.
Evidently, the Lion will have come within $\varepsilon/2$ of the Man
by time $t$ unless
%
%e31 ###
\begin{equation}\label{eqdiameter-bound}
\int^t_0 (1-\cos\alpha)\,\d s < \L(0) - \varepsilon/2
\leq
\diam_\mathrm{intr}(D) - \varepsilon/2
.
\end{equation}
Now consider the total curvature of the Lion's path. By Proposition~\ref
{propwell-posed}, the Lion's path
is uniformly approximated by pursuit paths driven by discretized
approximations to the Man's path.
If $x^{(n)}$ is the Lion's path driven by a discretized Man's path~$y^{(n)}$, then the Lion's path is
piecewise-geodesic, with total absolute curvature given by the sum of
the exterior angles formed
at the points that connect the geodesics that occur when $x^{(n)}$
changes direction. $\CAT0$ comparison
bounds then show the total curvature of $x^{(n)}$ is bounded above by
%
%e32 ###
\begin{equation}\label{eqTC-upper-bound}
\sum\frac{\sin\alpha^{(n)}}{\intrinsicdist(x^{(n)},y^{(n)})}\Delta
y^{(n)} ,
\end{equation}
where summation is over the jumps of the discretized path $y^{(n)}$,
and $\alpha^{(n)}$ is the exterior angle that
the jump $\Delta y^{(n)}$ contributes to the geodesic running from
$x^{(n)}$ to~$y^{(n)}$.

The total curvature of a path is a lower-semicontinuous function of the
path (using the uniform topology)
for $\CAT0$ spaces.
[This is a special case of a $\CAT\kappa$ result of \citet{KaruwannapatanaManeesawarng-2007}, referred to in
\citet{AlexanderBishopGhrist-2009}, Theorem 18.]
For the sake of completeness, we indicate the short proof for the $\CAT
0$ case. Consider a curve $q$
of finite length in a $\CAT0$ space. Its total curvature $\operatorname{TC}(q)$
is the supremum of sums of exterior angles
of piecewise-geodesic curves interpolating $q$; a $\CAT0$ comparison
argument shows that these sums of exterior angles
increase as the interpolating mesh is refined. Let $q^n$ be a sequence
of curves converging uniformly
to $q$. Furthermore, let $q^{n,m}$ be the piecewise-geodesic curve
interpolating $q^n$ at the points $k2^{-m}$ for $k=0, 1, \ldots.$
Then, by definition of total curvature,
\[
\operatorname{TC}(q^{n,m}) \nearrow \operatorname{TC}(q^n)\qquad \mbox{as }m\to\infty.
\]

\citeauthor{BridsonHaefliger-1999} [(\citeyear{BridsonHaefliger-1999}), Chapter II.3
Proposition~3.3]
observe that the $\CAT0$ property implies that \textit{interior}
angles are continuous functions of their end vertices and
upper-semicontinuous functions of their centre vertices. This
upper-semicontinuity
translates into lower-semicontinuity for exterior angles, and hence
\[
\limsup_{n\to\infty} \operatorname{TC}(q^{n,m}) \geq \operatorname{TC}(q^{\infty,m}) ,
\]
where $q^{\infty,m}$ is the uniform limit of $q^{n,m}$ as $n\to\infty$
[here we use the $\CAT0$ property again]
and is a piecewise-geodesic interpolation of $q$ at the points
$k2^{-m}$ for $k=0, 1, \ldots.$ Since
$\operatorname{TC}(q)=\lim\operatorname{TC}(q^{\infty,m})$, lower-semicontinuity now
follows from
\[
\limsup_{n\to\infty} \operatorname{TC}(q^{n}) \geq \limsup_{n\to\infty} \operatorname{TC}(q^{n,m})
\geq \operatorname{TC}(q^{\infty,m}) \to \operatorname{TC}(q)\qquad
\mbox{as } m\to\infty.
\]

Consequently, the upper bound \eqref{eqTC-upper-bound} provides an
upper bound on the total absolute curvature of the Lion's path in the limit.
Bearing in mind the $\operatorname{Lipschitz}(1)$ property of $y$, the total absolute
curvature $\tau(t)$ incurred by $x$
between times $0$ and~$t$ therefore satisfies
%
%e33 ###
\begin{equation}\label{eqtotal-curvature}
\tau(t)
\leq \int_0^t \frac{|\sin\alpha(s)|}{\L(s)}\,\d s .
\end{equation}
Assume that $\L(s) \geq\varepsilon/2$ for $s\leq t$.
By the Cauchy--Schwarz inequality
%KB and Inequality \eqref{eqdiameter-bound},
and~\eqref{eqdiameter-bound},
%
%e34 ###
\begin{eqnarray}\label{eqstep-on-the-way}
\tau(t)
&\leq&
\frac{2}{\varepsilon}\int_0^t|\sin\alpha|\,\d s
\leq
\frac{2}{\varepsilon}\sqrt{t \int_0^t \sin^2\alpha\,\d s}
\nonumber
\\[-8pt]
\\[-8pt]
\nonumber
&\leq&
\frac{2\sqrt{2}}{\varepsilon}\sqrt{t \int_0^t (1-\cos\alpha)\,\d s}
\leq
\frac{2\sqrt{2}}{\varepsilon}\sqrt{\diam_\mathrm{intr}(D)-\varepsilon
/2}\cdot\sqrt{t}.
\end{eqnarray}

Next, we follow \citet{AlexanderBishopGhrist-2006} in applying Reshetnyak
majorization [\citet{Reshetnyak-1968}; see also the telegraphic
description in
\citet{BerestovskijNikolaev-1993}, Section 7.4] to
generate a lower bound on the total absolute curvature of $\{x(s)\dvtx 0\leq
s\leq t\}$.
We provide details for the sake of completeness.

We argue as follows. Reshetnyak majorization asserts that for every
closed curve $\zeta$ in $\ol D$
[more generally, in any $\CAT0$ space],
one can construct a convex planar set $C$, bounded by a closed
unit-speed curve $\overline\zeta$,
and a distance-nonincreasing continuous map $\phi\dvtx C\to\ol D$ such
that $\phi\circ\overline\zeta=\zeta$;
moreover, $\phi$ preserves the arc-length distances along $\phi\circ
\overline\zeta$ and $\zeta$.
Consequently, $\phi$ restricted to $\partial C$ will not increase
angles and
the pre-images under $\phi$ of geodesic segments in $\zeta$ must
themselves be Euclidean geodesics
(i.e., line segments).

By our assumptions about $t$, the total absolute curvature of $\{
x(s)\dvtx 0\leq s\leq t\}$ is finite [see \eqref{eqstep-on-the-way}].
Fix an arbitrarily small $\delta_1\in(0,\pi/2)$.
It follows from the definitions of length and curvature of a path that,
% WSK3: Need some care here. Curves are parametrized by arc-length; this
% means the inscribed piecewise-geodesic curves $z$ will have different
%total lengths.
for each $n$,
we can approximate the unit-speed curve $\{x(s)\dvtx 0\leq s\leq t\}$
by a piecewise-geodesic curve $\{z(s)\dvtx 0\leq s\leq t'\}$
with the following properties:
\begin{itemize}
\item[--] The curve $z$ is parametrized using arc-length.
\item[--] Note that $x$ and $y$ are continuous, so is $\chi(x,y)$,
by Proposition~\ref{proplipschitz}(3). Hence, we can choose
$0=t_0<t_1<\cdots<t_n=t$ such that
the total absolute curvature of $\{x(s)\dvtx t_{i-1}\leq s\leq t_i\}$
is equal to $\pi/2 - \delta_1$ for all $i$, with the possible exception
of $i=n$.
\item[--] For every $i$, there exist $t_{i} = t_i^0 < t_i^1 < \cdots<
t_i^{m_i} = t_{i+1}$
and $s_i=s^{0}_i<s^{1}_i<\cdots<s_i^{m_i}=s_{i+1}$
such that $z(s_i^j) = x(t_i^j)$ and $z$ is geodesic on $[s_i^j,
s_i^{j+1}]$, for all $i$ and~$j$.
(Notice that the curve $z$ is inscribed in the curve $x$.)
\item[--] The total absolute curvature of $\{z(s)\dvtx s_{i-1}\leq s\leq
s_i\}$ is less than $\pi/2$.
In other words, the sum (over $j$) of exterior angles between $\{
z(s)\dvtx s_{i}^{j-1}\leq s\leq s_{i}^j\}$
and $\{z(s)\dvtx s_{i}^j\leq s\leq s_{i}^{j+1}\}$ at $s_i^j$ is less than
$\pi/2$.
[This is a consequence of $z$ being inscribed in $x$ and the $\CAT0$
property.]
\item[--] The difference between the lengths of $\{z(s)\dvtx 0\leq s\leq t'\}
$ and $\{x(s)\dvtx 0\leq s\leq t\}$ is less than $\delta_1$.
\end{itemize}
Then we have
%
%e35 ###
\begin{equation}\label{eqresh1}
\mbox{total absolute curvature}\bigl(\{x(s)\dvtx 0\leq s\leq t\}\bigr) \geq
\biggl(\frac{\pi}2 - \delta_1 \biggr) (n-1) .
\end{equation}

We apply Reshetnyak majorization to the closed curve formed by $\{
z(s)\dvtx\break s_{i-1}\leq s\leq s_i\}$
and its chord [the geodesic running from $z(t_i)$ back to $z(t_{i-1})$].
Reshetnyak majorization guarantees that the total absolute curvature of
$\{z(s)\dvtx s_{i-1}\leq s\leq s_i\}$
dominates the curvature of its pre-image in the boundary
of a convex planar set $C_i$.
Moreover, the perimeter of its pre-image in the boundary $C_i$ has length
$\length(\{z(s)\dvtx s_{i-1}\leq s\leq s_i\})$, while the remainder of the
boundary of $C_i$ must be
a line segment of length $\intrinsicdist(z(s_i), z(s_{i-1}))$.

The two-dimensional
pre-image of $\{z(s)\dvtx s_{i-1}\leq s\leq s_i\}$ therefore has total
curvature bound of $\frac{\pi}{2}$.
By two-dimensional Euclidean geometry, we can maximize
the ratio of the length of the pre-image of $\{z(s)\dvtx s_{i-1}\leq s\leq
s_i\}$
to the length $\intrinsicdist(z(s_i), z(s_{i-1}))$ of its chord
by considering the case of an isoceles right-angled triangle,
in which case the ratio is $\sqrt{2}$.
Accordingly, we obtain the upper bound
\[
\length\bigl(\{z(s)\dvtx s_{i-1}\leq s\leq s_i\}\bigr)
\leq
\sqrt{2} \intrinsicdist(z(s_i), z(s_{i-1}))
\leq
\sqrt{2} \diam_\mathrm{intr}(D) .
\]
It follows that a portion of the piecewise geodesic curve $z$ which
turns no more than $\frac{\pi}{2}$ cannot have length
exceeding $\sqrt{2}$ times the intrinsic diameter of the region. (Note
this is related to the Euclidean
diameter by Lemma~\ref{lemdiameter}.) This implies that we can control
the total length of $z$ and thus the total length of $x$,
with
%
%e36 ###
\begin{eqnarray}\label{eqresh2}
t-\delta_1 &=&
\length\bigl(\{x(s)\dvtx 0\leq s\leq t\}\bigr) -\delta_1
\nonumber
\\[-8pt]
\\[-8pt]
\nonumber
&\leq&
\length\bigl(\{z(s)\dvtx 0\leq s\leq t\}\bigr)
\leq
\sqrt{2} \diam_\mathrm{intr}(D) \times n .
\end{eqnarray}

Combining inequalities \eqref{eqresh1} and \eqref{eqresh2}, we deduce that
%
%e37 ###
\begin{eqnarray}\label{eqresh3}
&&\mbox{total absolute curvature}\bigl(\{
x(s)\dvtx 0\leq s\leq t\}\bigr)
\nonumber\\
&&\qquad\geq
\biggl(\frac{\pi}2 - \delta_1 \biggr)
(n-1)
\\
&&\qquad\geq
\biggl(\frac{\pi}2 - \delta_1 \biggr)
\biggl(\frac{t - \delta_1}
{\sqrt{2}\diam_\mathrm{intr}(D)}-1\biggr) .\nonumber
\end{eqnarray}
Recall that $\tau(t) = \mbox{total absolute curvature}(\{
x(s)\dvtx 0\leq s\leq t\})$ and $\length(\{x(s)\dvtx\break 0\leq s\leq t\})=t$. Letting
$\delta_1\to0$ in \eqref{eqresh3}, it follows that
\[
\frac{\tau(t)}{t} \geq
\frac{\pi}{2}\biggl(\frac{1}{\sqrt{2}\diam_\mathrm{intr}(D)}-\frac{1}{t}\biggr) .
\]
In combination with \eqref{eqstep-on-the-way}, this yields
\[
\frac{\pi}{2}\biggl(\frac{1}{\sqrt{2}\diam_\mathrm{intr}(D)}-\frac{1}{t}\biggr)
\leq
\frac{2\sqrt{2}}{\varepsilon}\sqrt{\diam_\mathrm{intr}(D)-\varepsilon
/2}\cdot\frac{1}{\sqrt{t}
}
\]
and hence the quadratic inequality for $q=\sqrt{t}$,
\[
\biggl(\frac{\pi}{2}\frac{1}{\sqrt{2}\diam_\mathrm{intr}(D)}\biggr) q^2
-
\biggl(\frac{2\sqrt{2}}{\varepsilon}\sqrt{\diam_\mathrm{intr}(D)-\varepsilon
/2}\biggr) q
-
\frac{\pi}{2}
\leq
0 .
\]

The left-hand side is negative for $q=0$ and the coefficient of $q^2$
is positive,
so there is exactly one positive root $q_c$
[which can be written out explicitly in terms of $\diam_\mathrm{intr}(D)$
and $\varepsilon$].
Combining this with our earlier arguments, it follows that
the Lion will come within $\varepsilon/2$ of the Man by time
$t_c := q_c^2$.~%
\end{pf*}

%s4 ###
\section{\texorpdfstring{From Brownian shy couplings to deterministic pursuit problems.}{From Brownian shy couplings to deterministic pursuit problems}}\label{seccoupling-to-pursuit}

This section is devoted to the proof of Theorem~\ref{thmno-shy-coupling}.
Consider a co-adapted coupling of reflecting Brownian motions $X$ and
$Y$ in
the bounded domain
$D\subseteq\Reals^d$
satisfying uniform exterior sphere and
interior cone conditions. \citet{Saisho-1987} showed that the reflected Brownian
motions can be realized by means of a Skorokhod transformation as
strong solutions
of stochastic differential equations driven by free Brownian motions.
As discussed in Section~\ref{secbasic-tools},
we can use
arguments embedded in the folklore of stochastic calculus, and stated explicitly
in \citet{Emery-2005} and in \citeauthor{Kendall-2009a} [(\citeyear{Kendall-2009a}), Lemma 6], to represent
this coupling as
%
%e39 ###
%e38 ###
\begin{eqnarray}
\d X &=& \d B - \nu_X\,\d L^X ,\label{j221}\\
\d Y &= &(\mathbb{J}^\top\,\d B+\mathbb{K}^\top\,\d A)-\nu_Y\,\d L^Y ,
\label{eqcoadaptedBM}
\end{eqnarray}
where $A$ and $B$ are independent $d$-dimensional Brownian motions, and
$\mathbb{J}$, $\mathbb{K}$ are predictable $(d\times d)$-matrix processes
such that
%
%e40 ###
\begin{equation}\label{m261}
\mathbb{J}^\top\mathbb{J}+\mathbb{K}^\top\mathbb{K} = (d\times d)\mbox{ identity matrix} .
\end{equation}
Here $L^X$ and $L^Y$ are the local times of $X$ and $Y$
on the boundary.\vadjust{\goodbreak}

The advantage of this explicit representation of the coupling is that
we can track
what happens to $X$ and $Y$ when we modify the Brownian motion $B$ by
adding a drift.
We will see that the effect of adding a very heavy drift based on the
vector field $\chi(X, Y)$ will be to convert
% Added equation reference following suggestion of referee in final
%remarks 04/10/2011
\eqref{j221} and \eqref{eqcoadaptedBM}
into a stochastic approximation of the deterministic Lion and Man
pursuit--evasion equations \eqref{eqdeterministic-skorokhod}
over a short time-scale.

%KB major changes starting here; for old version see previous file

\begin{prop}\label{o264} Suppose that $D\subset\Reals^d$ is $\CAT{0}$,
is bounded in the Euclidean metric,
and satisfies a uniform exterior sphere
condition and uniform interior cone condition.
For any $\eps>0$ and $X$ and $Y$ satisfying \eqref{j221} and \eqref{eqcoadaptedBM}
with $X(0)$, $Y(0)\in\overline{D}$, there exists $t>0$ such that
%
%e41 ###
\begin{equation}%\label{eqnonuniform-shyness}
\label{o264d}
\mathbb{P}\Bigl[\sup_{t/2\leq s\leq t} \idist(X(s),Y(s)) \leq\eps\Bigr] > 0 .
\end{equation}
\end{prop}
%
% WSK5: changed inf to sup as requested by MB.

\begin{pf}
Consider the following modification of \eqref{j221} and \eqref{eqcoadaptedBM},
%
%e43 ###
%e42 ###
\begin{eqnarray}
\label{o2620}X^n(t) &=& X(0)+B(t)
+\int_0^t n\chi(X^n(s), Y^n(s))\,\d s
\nonumber
\\[-9pt]
\\[-9pt]
\nonumber
&&{} - \int_0^t \nu_{X^n(s)}\,\d L^{X^n}_s
,\\[-2pt]
\label{o2621}Y^n(t) &=& Y(0) + \int_0^t \bigl(\mathbb{J}_s^\top\,\d B(s)+\mathbb{K}_s^\top
\,\d A(s)\bigr)
\nonumber
\\[-9pt]
\\[-9pt]
\nonumber
&&{}+ \int_0^t n
\mathbb{J}_s^\top\chi(X^n(s),Y^n(s)) \,\d s- \int_0^t \nu_{Y^n(s)}
\,\d L^{Y^n}_s .
\end{eqnarray}
By the Cameron--Martin--Girsanov theorem, the distributions of the
solutions of~\eqref{j221} and \eqref{eqcoadaptedBM} and \eqref{o2620} and \eqref{o2621} are
mutually absolutely continuous on every fixed finite interval.
We will show below that, after rescaling time, paths of $(X^n(\cdot),
Y^n(\cdot))$,
for large $n$, will be uniformly close to those for the corresponding
Lion and Man problem. Application
of Proposition~\ref{propgreedy} will then enable us to finish the proof.

We will make the following substitutions,
$X^n(t)=\widetilde{X}^n(nt)$,
$Y^n(t)=\widetilde{Y}^n(nt)$, $B(t)=\widetilde{B}^n(nt)/\sqrt{n}$,
$A(t)=\widetilde{A}^n(nt)/\sqrt{n}$,
$\mathbb{J}(t)={\widetilde{\mathbb{J}}^{(n)}}(nt)$,
$\mathbb{K}(t)={\widetilde{\mathbb{K}}^{(n)}}(nt)$.
Then \eqref{o2620} and \eqref{o2621} take the form
%
%e45 ###
%e44 ###
\begin{eqnarray}
\qquad\wt X^n(t) &=& X(0)+\frac{1}{\sqrt{n}} \wt B^n(t)
+\int_0^t \chi(\wt X^n(s),\wt Y^n(s))\,\d s - \int_0^t \nu_{\wt X^n(s)}\,\d L^{\wt X^n}_s ,\label{o2622}\\[-2pt]
\wt Y^n (t)&=& Y(0)+ \frac{1}{\sqrt{n}}
\int_0^t \bigl(\bigl(\widetilde{\mathbb{J}}_s^{(n)}\bigr)^\top\,\d
\widetilde{B}^n(s)+\bigl(\widetilde{\mathbb{K}}_s^{(n)}\bigr)^\top\,\d\widetilde
{A}^n(s)\bigr)\label{o2623}
\nonumber
\\[-9pt]
\\[-9pt]
\nonumber
&&{} +
\int_0^t\bigl(\widetilde{\mathbb{J}}_s^{(n)}\bigr)^\top
\chi(\widetilde{X}^n(s),
\widetilde{Y}^n(s))\,\d s- \int_0^t\nu_{\widetilde{Y}^n(s)}
\,\d L^{\widetilde{Y}^n_s} .
\end{eqnarray}
Note that $\wt B^n$ and $\wt A^n$ are Brownian motions.\vadjust{\goodbreak}

Consider the analog of (\ref{o2622}) and (\ref{o2623}), but without boundary:
%
%e47 ###
%e46 ###
\begin{eqnarray}
\wt U^n(t) &=& \frac{1}{\sqrt{n}} \wt B^n(t)
+\int_0^t \chi(\wt X^n(s),\wt Y^n(s))\,\d s ,\label{o2630}
\\
\label{o2631}\wt V^n (t)&=& \frac{1}{\sqrt{n}}
\int_0^t \bigl(\bigl(\widetilde{\mathbb{J}}_s^{(n)}\bigr)^\top\,\d
\widetilde{B}^n(s)+\bigl(\widetilde{\mathbb{K}}_s^{(n)}\bigr)^\top\,\d\widetilde
{A}^n(s)\bigr)
\nonumber
\\[-8pt]
\\[-8pt]
\nonumber
&&{}+ \int_0^t\bigl(\widetilde{\mathbb{J}}_s^{(n)}\bigr)^\top
\chi(\widetilde{X}^n(s),
\widetilde{Y}^n(s))\,\d s .
\end{eqnarray}
All components of the sextuplet
%
%e49 ###
%e48 ###
\begin{eqnarray}\label{o273}
\bX^n(t) &=& \biggl(\wt U^n(t) , \frac{1}{\sqrt{n}} \wt B^n(t) ,
\int_0^t \chi(\wt X^n(s),\wt Y^n(s))\,\d s ,
\nonumber\\
&&\phantom{\biggl(}{}\wt V^n (t), \frac{1}{\sqrt{n}}
\int_0^t \bigl(\bigl(\widetilde{\mathbb{J}}_s^{(n)}\bigr)^\top\,\d
\widetilde{B}^n(s)+\bigl(\widetilde{\mathbb{K}}_s^{(n)}\bigr)^\top\,\d\widetilde
{A}^n(s)\bigr) ,\\
&& \hspace*{95pt}{}\int_0^t\bigl(\widetilde{\mathbb{J}}_s^{(n)}\bigr)^\top
\chi(\widetilde{X}^n(s),
\widetilde{Y}^n(s))\,\d s \biggr)\nonumber
\end{eqnarray}
are tight by the criterion given by \citeauthor{StroockVaradhan-1979} [(\citeyear{StroockVaradhan-1979}), Section~1.4] since
the diffusion coefficients and the drifts are bounded by $1$. So, on an
appropriate
subsequence, $\bX^n$ converges weakly to a
limiting process $\bX^\infty$. By abuse of notation, we will denote
this subsequence $\bX^n$.
In particular, $\wt U^n(t)$ and $\wt V^n(t)$ converge
weakly, so, by
\citeauthor{Saisho-1987} [(\citeyear{Saisho-1987}), Theorem~4.1]
% WSK 8 July 2011 next line added to indicate why we need the
%conditions on $D$
(which applies because of the conditions imposed on $D$),
% WSK end
$(\wt X^n,\wt Y^n)$ converges weakly to a limiting continuous process
$(\wt X^\infty, \wt Y^\infty)$ along the same subsequence. It follows
that
%
%e50 ###
\begin{eqnarray}\label{o274}
\bY^n(t) &=& \biggl(\wt X^n(t),\wt Y^n(t), \wt U^n(t) , \frac{1}{\sqrt{n}} \wt
B^n(t) ,
\int_0^t \chi(\wt X^n(s),\wt Y^n(s))\,\d s ,
\nonumber\\
&&\phantom{\biggl(}{}\wt V^n (t), \frac{1}{\sqrt{n}}
\int_0^t \bigl(\bigl(\widetilde{\mathbb{J}}_s^{(n)}\bigr)^\top\,\d
\widetilde{B}^n(s)+\bigl(\widetilde{\mathbb{K}}_s^{(n)}\bigr)^\top\,\d\widetilde
{A}^n(s)\bigr) ,\\
&&\hspace*{118pt}{}\int_0^t\bigl(\widetilde{\mathbb{J}}_s^{(n)}\bigr)^\top
\chi(\widetilde{X}^n(s),
\widetilde{Y}^n(s))\,\d s \biggr)\nonumber
\end{eqnarray}
is tight and, therefore, converges weakly along a subsequence. Once again,
we will abuse the notation and assume that $\bY^n$ converges weakly. By the
Skorokhod lemma
%KB new
[\citet{EK}, Section~3.1, Theorem~1.8]
%KB end new
we can assume that the sequence $\bY^n$ converges
a.s., uniformly on compact intervals.
% MB1 - The notation for bX^n and bY^n is not so good, because it has
%no relationship with that
% for X^n and Y^n, and the internal relationship between bX^n and bY^n
%wrt is in fact different.
% Possibly use bold W and bold Z?
%KB bX and bY changed to bold W and Z
%I know what the Skorokhod Lemma is, but have not heard it
% referred to that way before; that may be the way it is now done, for
%all I know. Is this the standard
% designation?

The fourth and seventh components of $\bY^n$ are Brownian motions run
at rate
$\frac{1}{n}$ so they converge to the zero
process as $n \rightarrow\infty$.
The fifth and eighth components of $\bY^n$
are both $\operatorname{Lip}(1)$;
their limits are therefore also $\operatorname{Lip}(1)$. These observations
and \eqref{o2630} and \eqref{o2631} imply that the limits
$\wt V^\infty$ and $\wt U^\infty$ of $\wt V^n$ and $\wt U^n$ are $\operatorname{Lip}(1)$.

Let $\wt T^* = \inf\{t\geq0\dvtx  \wt X^\infty(t) = \wt Y^\infty(t)\}$.
The bounded vector field $\chi(\wt X^n, \wt Y^n)$ depends continuously
on $\wt X^n$ and $\wt Y^n$ (Proposition~\ref{proplipschitz}).
We may therefore apply the dominated convergence theorem
and \eqref{o2630} to deduce the following integral
representation for $\wt U^\infty$,
%
%e51 ###
\begin{equation}\label{o2660}
\wt U^\infty(t) = \int_0^t \chi(\wt X^\infty(s),\wt Y^\infty(s))\,\d s
\qquad\mbox{for } t < \wt T^*.
\end{equation}

Recall that, by the Skorokhod representation,
we can assume that $\wt X^n(t)$ and $\wt Y^n(t)$ converge almost surely.
Lemma~\ref{lemlipshitz} proved below shows
that $\wt X^\infty, \wt Y^\infty$
are both still $\operatorname{Lip}(1)$, \textit{with respect to the intrinsic metric
of $D$}.
Hence we can apply the results on $\CAT0$ Lion and Man
problems at the end of Section~\ref{secdeterministic-pursuit}.

Fix an arbitrarily small $\eps>0$. It follows from \eqref{o2660}
and from Proposition~\ref{propgreedy} that there exists $t_1 < \infty$
not depending on $X(0),Y(0)$ or $\omega$, such that $\idist(\wt X^\infty
(t),\wt Y^\infty(t)) \leq\eps/2$ for $t\geq t_1$.
We conclude that for some $n_0<\infty$, depending on $X(0)$ and $Y(0)$,
and all $n\geq n_0$,
\[%\label{o282}
\mathbb{P}\Bigl[
\sup_{t_1 \leq t \leq2t_1} \idist(\wt X^n(t),\wt Y^n(t)) \leq\eps
\Bigr] >0.
\]
Changing the clock to the original pace, we obtain
\[%\label{o282}
\mathbb{P}\Bigl[
\sup_{t_1/n \leq t \leq2t_1/n} \idist( X^n(t), Y^n(t)) \leq\eps
\Bigr] >0.
\]
%
%KB By the Girsanov theorem,
By the Cameron--Martin--Girsanov theorem,
% MB1 - Two pages back, the reference was to Cameron-Martin-Girsanovi,
%so presumably change this here.
%
%e52 ###
\begin{equation}\label{o282}
\mathbb{P}\Bigl[
\sup_{t_1/n \leq t \leq2t_1/n} \idist( X(t), Y(t)) \leq\eps
\Bigr] >0.
\end{equation}
Since $\eps>0$ is arbitrary, this completes the proof.
\end{pf}

%KB end of major changes; for old version see previous file

\begin{lem}\label{lemlipshitz}
Let $D$ be a domain satisfying uniform exterior sphere and interior cone
conditions.
Suppose that $Z$ is a continuous process on $\overline{D}$ derived by the
Skorokhod transformation from a free process $S$ that has $\operatorname{Lip}(1)$ sample
paths.
Then $Z$ itself has $\operatorname{Lip}(1)$ sample paths \textit{with respect to the
intrinsic metric}.
\end{lem}

\begin{pf}
Following \citeauthor{Saisho-1987} [(\citeyear{Saisho-1987}), Section~3], consider the step function
$S_m$ obtained from $S$ by sampling at instants $k2^{-m}$,
for $k=0, 1, \ldots.$ Suppose that $2^{-m}<r$, where
$r$ is the radius on which the uniform exterior sphere
condition is based.
Let $\overline{z}$ be the projection onto $\overline{D}$
described in Lemma~\ref{lemsg-for-uesc}.
The Skorokhod transformation\vadjust{\goodbreak} of $S_m$ is $Z_m$, given by
projecting increments back onto $\overline{D}$:
%
%e53 ###
\begin{equation}
Z_m(t) =
\cases{
\overline{Z_m\bigl((k-1)2^{-m}\bigr)+\Delta S_m(k2^{-m})}, \vspace*{2pt}\cr
\quad \hspace*{27pt}\mbox{for }
k2^{-m}\leq
t<(k+1)2^{-m},\vspace*{2pt}\cr
Z(0), \quad \mbox{for } 0\leq t<2^{-m}.}
\end{equation}
On account of the $\operatorname{Lip}(1)$ property of $S$, this
projection is defined when $2^{-m}<r$.

From \citeauthor{Saisho-1987} [(\citeyear{Saisho-1987}), Theorem 4.1], we know that $Z_m\to Z$
uniformly on
bounded time intervals. We compute the maximum possible Euclidean distance
between $Z_m(s)$ and $Z_m(t)$, if $0\leq t-s<2^{-m}$, when one or both of
$2^m s$, $2^m t$ are nonnegative integers.
Since $Z_m$ is constant on intervals $[k2^{-m},(k+1)2^{-m})$, it suffices
to produce an argument for the case when $2^m s=k-1$ and $2^mt=k$.
We therefore proceed to bound the Euclidean distance
$|Z_m(k2^{-m})-Z_m((k-1)2^{-m})|$. We will show that this can only exceed
$2^{-m}$ by an amount which, for large $m$, will make a negligible
contribution to path length when summed over the whole path.

If $Z_m(k2^{-m})\notin\partial D$, then there is nothing to prove,
since the
jump is $\Delta S_m(k2^{-m})$,
which
is bounded in length by $2^{-m}$
since $S$ is $\operatorname{Lip}(1)$. So we instead suppose that $Z_m(k2^{-m})\in
\partial D$.
For convenience, set $y=Z_m(k2^{-m})-(Z_m((k-1)2^{-m})+\Delta
S_m(k2^{-m}))$ to
be the Skorokhod correction to be applied at this step, and set
$a=|\Delta
S_m(k2^{-m})|$ to be the length of the uncorrected jump. Finally, let
$\theta$ be the angle between the vector $y$ and the negative jump
$-\Delta
S_m(k2^{-m})$. These definitions are illustrated in Figure~\ref{figshy2-max},
together with the supporting ball~$B$ at $Z_m(k2^{-m})\in\partial D$ whose
centre is located at $Z_m(k2^{-m})-\lambda y$ for some $\lambda
=r/|y|>0$ and whose
existence is guaranteed by the construction of the $x\mapsto\overline{x}$
projection map as described in Lemma~\ref{lemsg-for-uesc}.

%f9 ###
\begin{figure}[b]

\includegraphics{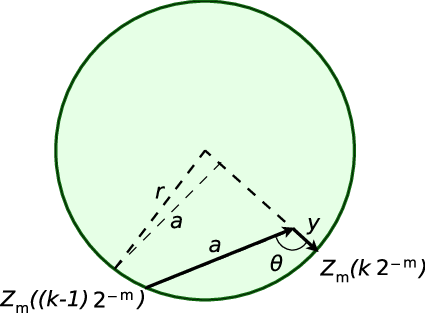}

\caption{Illustration of the geometry underlying the argument of Lemma
\protect\ref{lemlipshitz}.}\label{figshy2-max}
\end{figure}

First note that
$|Z_m(k2^{-m})-Z_m((k-1)2^{-m})|=\sqrt{a^2+y^2-2ay\cos\theta}$\break
(where we abuse notation by letting $y$ also stand for the length of
the vector $y$). This
increases as $\theta$ increases to $\pi$, so long as $a$, $y$,
$Z_m(k2^{-m})$ are held fixed. Thus we can assume that $\theta$ has increased
to the point where $Z_m((k-1)2^{-m})$, as well as $Z_m(k2^{-m})$,
belong to
$\partial B$. (This will happen if, as required above, $2^{-m}<r$.) Now
observe that the distance $|Z_m(k2^{-m})-Z_m((k-1)2^{-m})|$ will be bounded
above by the smaller of the two distances from $Z_m(k2^{-m})$ to the
intercepts of $\partial B$ by a line parallel to $y$, and at distance $a$
from $Z_m(k2^{-m})$. But two applications of Pythagoras' theorem show that
this distance is given by
\begin{eqnarray*}
\sqrt{a^2 +\bigl(r-\sqrt{r^2-a^2}\bigr)^2} &=& \sqrt{2r^2-2r\sqrt{r^2-a^2}}\\
&=& \sqrt{2}r\sqrt{1-\sqrt{1-\frac{a^2}{r^2}}}
=
\sqrt{2}r\sqrt{\frac{a^2}{2r^2}\biggl(1+\frac{1}{4}(z^*)^2\biggr)}
\end{eqnarray*}
for some $z^*$ in the range $[0,\frac{r^2}{a^2}]$. (The last step arises
from a second-order Taylor series expansion.) Therefore
\[
\sqrt{a^2 +\bigl(r-\sqrt{r^2-a^2}\bigr)^2} \leq
\sqrt{2}r\sqrt{\frac{a^2}{2r^2}\biggl(1+\frac{a^2}{4r^2}\biggr)}
\leq
a \biggl(1+\frac{a^2}{8r^2}\biggr)
\]
(using $\sqrt{1+z}\leq1+\frac{1}{2}z$ for $z\geq-1$).

Thus the total path length over the time interval $(s,t)$ is bounded
above by
\[
\bigl((t-s)2^m+2\bigr)\times2^{-m}\biggl(1+\frac{2^{-2m}}{8r^2}\biggr),
\]
which converges to $t-s$ as $m\to\infty$.
Hence we obtain
\[
\intrinsicdist(Z(s),Z(t))
\leq
t-s ,
\]
thus establishing the $\operatorname{Lip}(1)$ property in intrinsic metric for $Z$.
\end{pf}

%KB more material deleted
%Proposition~\ref{propgreedy}
%states that $X^\infty$ comes within distance %$\varepsilon/2$ of
%$Y^\infty$ in the intrinsic metric by time %$\tau=t_c$. Since
%$(\widetilde{X}^n, \widetilde{Y}^n)$ converges weakly %to $(X^\infty,
%Y^\infty)$,
%it follows that, for sufficiently large
%$n$, there is a positive probability of %$\widetilde{X}^n$ coming
%within
%distance $\varepsilon$ of $\widetilde{Y}^n$ by time %$\tau=t_c$ in
%this metric. On the other hand,
%$\intrinsicdist(x,y)\geq\|x-y\|$. Therefore, by the
%%Cameron--Martin--Girsanov
%argument given at the beginning of the section, for $X$ %and $Y$
%satisfying \eqref{eqcoadaptedBM}
%and given $X(0)$, $Y(0)\in\overline{D}$, there is a %positive
%probability of $X$ coming within distance %$\varepsilon$ of $Y$ by
%time $t=t_c/n$.
We will show that the bound in Proposition~\ref{o264} is uniform over
all $X(0)$ and $Y(0)$. We will switch from the intrinsic distance to
the Euclidean distance in the formulation of the next proposition. This
is legitimate in view of
\eqref{m232}.

% WSK3: delete
% (which we state in a slightly more general form so that it can be
%re-used
% in Section~\ref{secstar-shaped}).

\begin{prop}\label{m263}
% WSK5 Next line added as requested by Maury.
Let $D$ be a domain satisfying uniform exterior sphere and interior
cone conditions.
% Let $K\subseteq\overline{D}$ be a closed set.
Suppose that there exist $t_1>0$ and $\varepsilon_1>0$ such that,
for any $X$ and $Y$ satisfying
% Added equation reference following suggestion of referee in final
%remarks 04/10/2011
\eqref{j221} and \eqref{eqcoadaptedBM}
with $X(0)$, $Y(0)\in\overline{D}$,
%
%e54 ###
\begin{equation}\label{eqnonuniform-shyness}
\mathbb{P}\Bigl[\inf_{0\leq t\leq t_1} |X(t)-Y(t)| \leq\varepsilon_1\Bigr] > 0 .
\end{equation}
Then
%
%e55 ###
\begin{equation}\label{equniform-shyness}
\mathbb{P}\Bigl[\inf_{0\leq t\leq t_1} |X(t)-Y(t)| \leq\varepsilon_1\Bigr] > p_1
\end{equation}
for some $p_1>0$ not depending on $X(0)$ and $Y(0)$.\vadjust{\goodbreak}
\end{prop}

\begin{pf}
Suppose \eqref{equniform-shyness} does not hold. Then there exist
$t_1>0$, $\varepsilon_1>0$, sequences $\{x_n\}_{n\geq1}$,
$\{y_n\}_{n\geq1}$ of points in $\overline{D}$, random processes $\{
A_t, t\geq0\}$,
$\{B_t, t\geq0\}$, $\{\JJ^n_t, t\geq0\}$ and $\{\KK^n_t, t\geq0\}$,
and solutions $X^n$ and $Y^n$ of
% Added equation reference following suggestion of referee in final
%remarks 04/10/2011
\eqref{j221} and \eqref{eqcoadaptedBM}
satisfying the following properties.
The processes $A$ and $B$ are $d$-dimensional Brownian motions starting
from 0, and independent of each other.
The $(d\times d)$-matrix-valued processes $\JJ^n$ and $\KK^n$ are
predictable with respect to the natural filtration of $A$ and~$B$,
such that $(\JJ^n_t)^\top\JJ_t^n+(\KK^n_t)^\top\KK_t$ is the $(d\times
d)$ identity matrix at all times $t$.
Let $X^n$ and $Y^n$ denote solutions to
% Added equation reference following suggestion of referee in final
%remarks 04/10/2011
\eqref{j221} and \eqref{eqcoadaptedBM}
based on the Brownian motions $A$ and~$B$, using the predictable
integrators $\JJ^n$ and $\KK^n$,
and starting from $X^n(0)=x_n\in\overline{D}$ and $Y^n(0)=\break y_n\in
\overline{D}$. Then
%
%e56 ###
\begin{equation}\label{m262}
\mathbb{P}\Bigl[\inf_{0\leq t \leq t_1} |X^n(t) - Y^n(t)| > \eps_1 \Bigr] > 1-2^{-n}.
\end{equation}

Let
$(M^{n,1}_t, M^{n,2}_t) = (\int_0^t \,\d B_s, \int_0^t (\JJ^n_s)^\top\,\d B_s+ \int_0^t (\KK^n_s)^\top\,\d A_s)$.
The processes $M^{n,1}$ and $M^{n,2}$ are Brownian motions and so the
sequence of pairs is tight, which therefore
possesses a subsequence converging in distribution. By abuse of
notation, we assume that
the whole sequence $(M^{n,1}, M^{n,2})$ converges in
distribution to, say, $(M^{\infty,1}, M^{\infty,2})$. It is clear that
$M^{\infty,1}$ and $ M^{\infty,2}$
are Brownian motions.

Let $\FF_t = \sigma((M^{\infty,1}_s, M^{\infty,2}_s), s\leq t)$ be the
natural filtration for
$(M^{\infty,1}, M^{\infty,2})$. We will show that $(M^{\infty,1},
M^{\infty,2})$ are co-adapted Brownian motions relative to $\{\FF_t\}$.
Since $(M^{n,1}, M^{n,2})$ are co-adapted Brownian motions, for all
$0\leq t_1 \leq t_2 \leq\cdots\leq t_n \leq t \leq s_1 \leq s_2$,
the random variable $M^{n,1}_{s_2} - M^{n,1}_{s_1}$
is independent of
\[
((M^{n,1}_{t_1}, M^{n,2}_{t_1}), (M^{n,1}_{t_2}, M^{n,2}_{t_2}), \ldots,
(M^{n,1}_{t_n}, M^{n,2}_{t_n})).
\]
Independence is preserved by weak limits, so
$M^{\infty,1}_{s_2} - M^{\infty,1}_{s_1}$
is independent of
\[
((M^{\infty,1}_{t_1}, M^{\infty,2}_{t_1}), (M^{\infty,1}_{t_2},
M^{\infty,2}_{t_2}), \ldots, (M^{\infty,1}_{t_n}, M^{\infty,2}_{t_n})).
\]
This implies that $M^{\infty,1}_{s_2} - M^{\infty,1}_{s_1}$ is
independent of $\FF_t$. Since the same argument applies to $M^{\infty
,2}_{s_2} - M^{\infty,2}_{s_1}$, we see that $(M^{\infty,1}, M^{\infty
,2})$ are co-adapted relative to $\{\FF_t\}$.
Recall from Section~\ref{secbasic-tools} that this implies that there
exist Brownian motions $\{A^\infty_t, t\geq0\}$ and $\{B^\infty_t,
t\geq0\}$ and
processes $\{\JJ^\infty_t, t\geq0\}$ and $\{\KK^\infty_t, t\geq0\}$
such that $(M^{\infty,1}_t, M^{\infty,2}_t) = (\int_0^t \,\d B^\infty_s,
\int_0^t (\JJ^\infty_s)^\top\,\d B^\infty_s+ \int_0^t (\KK^\infty_s)^\top
\,\d A^\infty_s)$.

Recall that $(M^{n,1}, M^{n,2}) \to(M^{\infty,1}, M^{\infty,2})$
weakly in the uniform topology on all compact intervals. Going back to
the original notation, we see that
\begin{eqnarray*}
&&\biggl(\int_0^t \,\d B_s, \int_0^t (\JJ^n_s)^\top\,\d B_s+ \int_0^t (\KK
^n_s)^\top\,\d A_s\biggr)\\
&&\qquad \to\biggl(\int_0^t \,\d B^\infty_s, \int_0^t (\JJ^\infty_s)^\top\,\d B^\infty_s+
\int_0^t (\KK^\infty_s)^\top\,\d A^\infty_s\biggr)
\end{eqnarray*}
weakly in the uniform topology on all compact intervals.
By the Skorokhod lemma, we can assume that the processes converge a.s.
in the supremum topology on compact intervals.

% WSK5 added clause about continuous dependence below
Since
$\ol D$ is compact, we can assume, passing to a subsequence if
necessary, that the initial points satisfy
$x_n \to x_\infty\in\overline{D}$ and $y_n \to y_\infty\in\ol D$ as
$n\to\infty$.
In view of the representation of coupled reflected Brownian motions
using stochastic differential equations
% Added equation reference following suggestion of referee in final
%remarks 04/10/2011
\eqref{j221} and \eqref{eqcoadaptedBM},
established in \citeauthor{Saisho-1987} [(\citeyear{Saisho-1987}), Theorem 4.1], and employing the
continuous dependence on driving Brownian motions established there, we
see that
$(X^n, Y^n) \to(X^\infty, Y^\infty)$ weakly in the uniform topology on
all compact intervals,
where $(X^\infty, Y^\infty)$ represents the solution to
% Added equation reference following suggestion of referee in final
%remarks 04/10/2011
\eqref{j221} and \eqref{eqcoadaptedBM}
with $X^\infty(0) = x_\infty$,
$Y^\infty(0) = y_\infty$, corresponding to $A^\infty$, $B^\infty$, $\JJ
^\infty$ and $\KK^\infty$.
We obtain from \eqref{m262} and weak convergence of
$(X^n, Y^n)$ to $(X^\infty, Y^\infty)$ that, for every $n$,
\[%\label{m262}
\mathbb{P}\Bigl[\inf_{0\leq t \leq t_1} |X^\infty(t) - Y^\infty(t)| \geq\eps_1
\Bigr]\geq 1 - 2^{-n}.
\]
Taking the limit as $n\to\infty$, this contradicts \eqref
{eqnonuniform-shyness} in the statement of the Proposition.
Consequently \eqref{equniform-shyness} must hold for some $p_1$.
\end{pf}

We now complete the proof of Theorem~\ref{thmno-shy-coupling}, applying
Proposition~\ref{m263}
% WSK3: delete
% (in the special case of $K=\ol D$)
together with standard reasoning.
Consider processes $X$ and $Y$ starting from any pair of points in $\ol
D$ and corresponding to any ``strategy'' $\JJ$ and $\KK$.
Because of the uniform bound in Proposition~\ref{m263}, the probability
of $X$ and $Y$ not coming within
distance~$\eps_1$ of each other on the interval $[kt_1, (k+1)t_1]$,
conditional on not coming within this distance
before $kt_1$, is bounded above by $1-p_1$ for any $k$, by the Markov
property. Hence, the probability
of $X$ and $Y$ not coming within distance $\eps_1$ of each other on the
interval $[0, kt_1]$ is bounded above
by $(1-p_1)^k$. Letting $k\to\infty$, it follows that $X$ and $Y$ are
not $\eps_1$-shy. Since $\eps_1$ can be
taken arbitrarily small, the proof of Theorem~\ref{thmno-shy-coupling}
is complete.

We remark that the matrices $\mathbb{J}$ and $\mathbb{K}$ employed in
% Added equation reference following suggestion of referee in final
%remarks 04/10/2011
\eqref{j221} and \eqref{eqcoadaptedBM}
are predictable
and, consequently,
%the choice of the pursuer's velocity is made without
%looking into the future.
the choice of the pursuer's velocity is based strictly on past information.
This is in
contrast to the pursuit--evasion problems and associated paradoxes
discussed by \citet{BollobasLeaderWalters-2009}.

%s5 ###
\section{\texorpdfstring{Complements and conclusions.}{Complements and conclusions}}\label{seccomplements}
We conclude this paper by remarking on some supplementary results and concepts,
and by considering possibilities for future work.

%s5.1 ###
\subsection{\texorpdfstring{Comparison with previous methods.}{Comparison with previous methods}}
The fundamental idea in this paper turns out in the end to resemble
that of
\citet{BenjaminiBurdzyChen-2007}, but uses simple notions of weak convergence
and tightness, rather than detailed large deviation estimates.\vadjust{\goodbreak}
Moreover, the use
of metric geometry notions enables us to finesse many analytical technicalities.
(Perhaps this is the first application of modern metric geometry to Euclidean
stochastic calculus?) On the other hand, the stochastic control methods of
\citet{Kendall-2009a} are quite different. The stochastic control
approach uses
potential theory to estimate the value function of an associated stochastic
game; consequently the methods of \citet{Kendall-2009a} may be expected
to give
sharper information (bounds on expectation of stopping times), but in more
limited cases (convexity of domain). However, one can observe that, at
least in
principle, the stochastic game formulation still applies in the general case.
For example, there is a value function to be discovered for a stochastic
control reformulation of Theorem~\ref{thmno-shy-coupling}, and in
principle it
might be possible to estimate this value function and so gain more information
than is supplied by the weak geometric bounds established above.

We note that
many promising ideas based on stochastic calculus fail to show
nonshyness because they cannot be applied to ``perverse''
couplings with the property that, on some time intervals,
$|X-Y|$ grows at a deterministic rate [see Example 4.2 of
\citet{BenjaminiBurdzyChen-2007}].

Also note that the proof in \citet{Kendall-2009a}, which works in
convex domains, does
not appear to be (directly) extendable to calculations involving the intrinsic
metric---simple manipulation using symbolic It\^o calculus [\citet{Kendall-2001b}] shows that
the drift of
$\intrinsicdist(X,Y)$
is unbounded at distances bounded away from zero. In particular,
Bessel-like divergences
for $\intrinsicdist(X,Y)$ of magnitude $a$ occur when the
geodesic from $X$ to $Y$ touches a concave part of $\partial D$
at $x$ and $|x-Y|= 1/a$. The first-order differential geometry given in
Proposition~\ref{propgauss} (the generalization of Gauss' lemma)
is the best we can do for $\CAT0$ domains satisfying uniform exterior sphere
and interior cone conditions.

% \subsection[Higher dimensions and the failure of CAT0]{Higher
%dimensions and the failure of $\CAT0$}
%s5.2 ###
\subsection{\texorpdfstring{Higher dimensions and the
failure of $\CAT0$.}{Higher dimensions and the
failure of $\CAT0$}}

For planar domains, $\CAT0$ and simple-connectedness are equivalent,
in which case, by Theorem~\ref{thmno-shy-coupling2}, there are no shy
co-adapted couplings.
% WSK3: deleted
% As illustrated by Theorem~\ref{thmstar-shaped}, appropriate
%conditions for ruling
% out shy couplings are less clear in higher dimensions. Based on these
%considerations,
% it is natural to consider replacing the simply-connected condition by
%contractibility in
% general dimensions. Accordingly, we formulate a bold and possibly
%rash conjecture:
% WSK3: add (then WSK4 delete!)
% In higher dimensions, it is natural to ask whether the \CAT0
%condition is essential for there to be no shy coupling.
% We do not at all believe this to be the case. It is possible to give
%an argument that suggests
% that star-shaped domains with smooth boundary conditions cannot
%support shy coupling, by establishing the analogous
% result for a corresponding deterministic pursuit-evasion problem. To
%apply this argument to the probabilistic case would
% require careful potential-theoretic or large-deviation arguments; we
%therefore leave this as a project for another day.
% WSK4: add
In higher dimensions, it is natural to ask whether the $\CAT0$ condition
is essential for there to be no shy coupling.
We do not at all believe this to be the case. It is possible to give an
argument suggesting
that star-shaped domains with smooth boundary conditions cannot support
shy couplings, by establishing the analogous
result for a corresponding deterministic pursuit--evasion problem. To
apply this argument to the probabilistic case would
require more careful arguments. We therefore leave this as a project
for another day.

As a spur to future work, we formulate a bold and possibly rash conjecture:
\begin{conjecture}\label{conjcontractible}
There can be no shy co-adapted coupling for reflecting Brownian motions in
bounded contractible domains in any dimension.
\end{conjecture}
%
% WSK3: deleted
% It is certainly possible to combine the techniques of Theorems
% to deal with some non-$\CAT0$, nonstar-shaped domains on an
% WSK3: added

While resolution of the star-shaped case appears to be largely a
technical matter,
we believe that new ideas will be required to make substantial progress
toward resolving the conjecture.

%s5.3 ###
\subsection{When can shyness exist?}\label{secshyness-exist}
Many examples of shy couplings can be generated using suitable symmetries.
However, we do not know of any examples in which symmetries play no r\^ole.
Accordingly we formulate a further conjecture:
\begin{conjecture}\label{conjsymmetry}
If a bounded domain $D$ supports a shy co-adapted coupling for reflecting
Brownian motions, then there exists a shy co-adapted coupling that can
be realized using a
rigid-motion
symmetry of the
domain~$D$.
\end{conjecture}

%KB new
A stronger form of the above conjecture, saying ``If a bounded domain
$D$ supports a shy co-adapted coupling for reflecting
Brownian motions, then the shy coupling is realized using a
rigid-motion symmetry of the domain~$D$,'' is false. To see this,
consider the planar annulus $A = \ball(0,2) \setminus\ball(0,1)$ and
let $\TT$ be the symmetry with respect the origin. Let $X$ be reflected
Brownian motion in $A$ and $Y = \TT(X)$. Let $D = A \times(0,1)$ and
let $Z$ be reflected Brownian motion in $(0,1)$, independent of $X$ and
$Y$. Then $(X,Z)$ and $(Y,Z)$ form a shy coupling in $D$ which cannot
be realized using a rigid-motion symmetry of $D$.
%KB end new

Note
that \citeauthor{BenjaminiBurdzyChen-2007} [(\citeyear{BenjaminiBurdzyChen-2007}), Example 3.9] supplies an example
based on
Brownian motion on graphs, for which there is no fixed-point-free
isometry and yet
a shy coupling exists.
However we do not see how to use the idea of this construction
to construct a counterexample to the above conjecture.

%s5.4 ###
\subsection{\texorpdfstring{Further questions.}{Further questions}}

We enumerate a short list of additional questions.
\begin{longlist}[(1)]
\item[(1)] Shyness is interesting for foundational reasons: coupling is an important
tool in probability, and shyness informs us about coupling. We do not
know of
any honest applications of shyness. However, one can contrive a kind of
cryptographic context. Suppose one wishes to mimic a target $Y$, which
is a
randomly evolving high-dimensional structure, in such a way that the
mimic $X$
never comes within a certain distance of the target $Y$.
Shyness concerns the question, whether it is possible to do this in a
way that
is perfectly concealed from an observer watching the mimic $X$
alone.
\item[(2)] In this formulation, it is not clear why one should restrict
consideration to co-adapted
couplings. Our methods do not lend themselves to the non-co-adapted
case, and
the question is open whether or not results change substantially if one is
allowed to use such couplings. In particular, it seems possible that
Conjecture~\ref{conjsymmetry}
might have a quite different answer in this context.
%KB \item In further work we plan to develop these ideas in %the
%context of deterministic pursuit-evasion, with %implications for the
%absence
%of shy-coupling for star-shaped domains.
% \item In further work (\cite{BBK}) we plan to
% study deterministic pursuit-evasion problem and shy couplings
% in multidimensional star-shaped domains.
%
\item[(3)] In further work [\citet{BBK}] we plan to study the
deterministic pursuit--evasion problem,
in conjunction with shy couplings, for multidimensional $\CAT\kappa$
domains possessing ``stable rubber bands,'' a condition that is partly
topological and partly geometric. As a corollary, we plan to prove that
there are no shy couplings in multidimensional star-shaped domains.
\item[(4)] The Lion and Man problem has been generalized to the case of
multiple Lions. [An
early instance is given in \citet{Croft-1964}.] Can one formulate and prove
useful results for a corresponding notion of multiple shyness?
\end{longlist}
%

%KB \noindent{\bf Acknowledgments}. We are grateful to Soumik Pal for
%giving most helpful advice.
\section*{\texorpdfstring{Acknowledgments.}{Acknowledgments}}
We are grateful to Stephanie Alexander, Richard\break Bishop, Chanyoung Jun
and Soumik Pal for giving most helpful advice.
The presentation of our results was improved by suggestions from the
referee and associate editor for which we are thankful.

% imsref loaded by akundreckaite, 2012-03-28 12:47:35
%

\printaddresses

\end{document}